\input epsf
\input amstex
\input amsppt.sty

\magnification1200
\hsize14cm
\vsize19cm

\TagsOnRight

\def\ZeilAV{20}
\def\ZeilAM{19}
\def\StemAE{18}
\def\SlatAC{17}
\def\PrBMAA{16}
\def\PropAA{15}
\def\PeWZAA{14}
\def\PaScAA{13}
\def\MM{12}
\def\KrZeAA{11}
\def\KratBD{10}
\def\KratBI{9}
\def\KratBH{8}
\def\KratBG{7}
\def\KratAM{6}
\def\HeGeAA{5}
\def\GeViAB{4}
\def\GeViAA{3}
\def\DT{2}
\def\CiucAB{1}


\def\FZ{10}

\def\Aa{1.1}

\def\AC{1.2}
\def\Ac{2.3}
\def\AE{2.4}
\def\AF{2.5}
\def\BA{2.1}
\def\BB{2.2}

\def\DA{4.1}
\def\DAa{4.2}
\def\DD{4.3}

\def\EB{5.1}
\def\EC{5.2}
\def\EE{5.3}
\def\EF{5.4}
\def\EG{5.5}
\def\EH{5.6}
\def\EI{5.7}
\def\EIa{5.8}
\def\EJ{5.9}
\def\EKa{5.10}
\def\EKb{5.11}
\def\ELa{5.12}
\def\ELb{5.13}
\def\EMa{5.14}
\def\EMb{5.15}
\def\EN{5.16}
\def\EO{5.17}
\def\EPa{5.18}
\def\EPb{5.19}
\def\EQa{5.20}
\def\EQb{5.21}
\def\EIb{5.22}
\def\ERa{5.23}
\def\ERb{5.24}
\def\ES{5.25}
\def\ET{5.26}

\def\TA{1}
\def\TB{2}
\def\TC{3}
\def\TD{4}

\def\TK{11}
\def\TKa{12}
\def\TL{13}
\def\TM{14}
\def\TN{15}
\def\TO{16}

\def\po#1#2{(#1)_#2}
\def\({\left(}
\def\){\right)}

\catcode`\@=11
\font\tenln    = line10
\font\tenlnw   = linew10

\newskip\Einheit \Einheit=0.5cm
\newcount\xcoord \newcount\ycoord
\newdimen\xdim \newdimen\ydim \newdimen\PfadD@cke \newdimen\Pfadd@cke

\newcount\@tempcnta
\newcount\@tempcntb

\newdimen\@tempdima
\newdimen\@tempdimb

\newdimen\@wholewidth
\newdimen\@halfwidth

\newcount\@xarg
\newcount\@yarg
\newcount\@yyarg
\newbox\@linechar
\newbox\@tempboxa
\newdimen\@linelen
\newdimen\@clnwd
\newdimen\@clnht

\newif\if@negarg

\def\@whilenoop#1{}
\def\@whiledim#1\do #2{\ifdim #1\relax#2\@iwhiledim{#1\relax#2}\fi}
\def\@iwhiledim#1{\ifdim #1\let\@nextwhile=\@iwhiledim
        \else\let\@nextwhile=\@whilenoop\fi\@nextwhile{#1}}

\def\@whileswnoop#1\fi{}
\def\@whilesw#1\fi#2{#1#2\@iwhilesw{#1#2}\fi\fi}
\def\@iwhilesw#1\fi{#1\let\@nextwhile=\@iwhilesw
         \else\let\@nextwhile=\@whileswnoop\fi\@nextwhile{#1}\fi}

\def\thinlines{\let\@linefnt\tenln \let\@circlefnt\tencirc
  \@wholewidth\fontdimen8\tenln \@halfwidth .5\@wholewidth}
\def\thicklines{\let\@linefnt\tenlnw \let\@circlefnt\tencircw
  \@wholewidth\fontdimen8\tenlnw \@halfwidth .5\@wholewidth}
\thinlines

\PfadD@cke1pt \Pfadd@cke0.5pt
\def\PfadDicke#1{\PfadD@cke#1 \divide\PfadD@cke by2 \Pfadd@cke\PfadD@cke \multiply\PfadD@cke by2}
\long\def\LOOP#1\REPEAT{\def\BODY{#1}\ITERATE}
\def\ITERATE{\BODY \let\next\ITERATE \else\let\next\relax\fi \next}
\let\REPEAT=\fi
\def\Punkt{\hbox{\raise-2pt\hbox to0pt{\hss$\ssize\bullet$\hss}}}
\def\DuennPunkt(#1,#2){\unskip
  \raise#2 \Einheit\hbox to0pt{\hskip#1 \Einheit
          \raise-2.5pt\hbox to0pt{\hss$\bullet$\hss}\hss}}
\def\NormalPunkt(#1,#2){\unskip
  \raise#2 \Einheit\hbox to0pt{\hskip#1 \Einheit
          \raise-3pt\hbox to0pt{\hss\twelvepoint$\bullet$\hss}\hss}}
\def\DickPunkt(#1,#2){\unskip
  \raise#2 \Einheit\hbox to0pt{\hskip#1 \Einheit
          \raise-4pt\hbox to0pt{\hss\fourteenpoint$\bullet$\hss}\hss}}
\def\Kreis(#1,#2){\unskip
  \raise#2 \Einheit\hbox to0pt{\hskip#1 \Einheit
          \raise-4pt\hbox to0pt{\hss\fourteenpoint$\circ$\hss}\hss}}

\def\Line@(#1,#2)#3{\@xarg #1\relax \@yarg #2\relax
\@linelen=#3\Einheit
\ifnum\@xarg =0 \@vline
  \else \ifnum\@yarg =0 \@hline \else \@sline\fi
\fi}

\def\@sline{\ifnum\@xarg< 0 \@negargtrue \@xarg -\@xarg \@yyarg -\@yarg
  \else \@negargfalse \@yyarg \@yarg \fi
\ifnum \@yyarg >0 \@tempcnta\@yyarg \else \@tempcnta -\@yyarg \fi
\ifnum\@tempcnta>6 \@badlinearg\@tempcnta0 \fi
\ifnum\@xarg>6 \@badlinearg\@xarg 1 \fi
\setbox\@linechar\hbox{\@linefnt\@getlinechar(\@xarg,\@yyarg)}%
\ifnum \@yarg >0 \let\@upordown\raise \@clnht\z@
   \else\let\@upordown\lower \@clnht \ht\@linechar\fi
\@clnwd=\wd\@linechar
\if@negarg \hskip -\wd\@linechar \def\@tempa{\hskip -2\wd\@linechar}\else
     \let\@tempa\relax \fi
\@whiledim \@clnwd <\@linelen \do
  {\@upordown\@clnht\copy\@linechar
   \@tempa
   \advance\@clnht \ht\@linechar
   \advance\@clnwd \wd\@linechar}%
\advance\@clnht -\ht\@linechar
\advance\@clnwd -\wd\@linechar
\@tempdima\@linelen\advance\@tempdima -\@clnwd
\@tempdimb\@tempdima\advance\@tempdimb -\wd\@linechar
\if@negarg \hskip -\@tempdimb \else \hskip \@tempdimb \fi
\multiply\@tempdima \@m
\@tempcnta \@tempdima \@tempdima \wd\@linechar \divide\@tempcnta \@tempdima
\@tempdima \ht\@linechar \multiply\@tempdima \@tempcnta
\divide\@tempdima \@m
\advance\@clnht \@tempdima
\ifdim \@linelen <\wd\@linechar
   \hskip \wd\@linechar
  \else\@upordown\@clnht\copy\@linechar\fi}

\def\@hline{\ifnum \@xarg <0 \hskip -\@linelen \fi
\vrule height\Pfadd@cke width \@linelen depth\Pfadd@cke
\ifnum \@xarg <0 \hskip -\@linelen \fi}

\def\@getlinechar(#1,#2){\@tempcnta#1\relax\multiply\@tempcnta 8
\advance\@tempcnta -9 \ifnum #2>0 \advance\@tempcnta #2\relax\else
\advance\@tempcnta -#2\relax\advance\@tempcnta 64 \fi
\char\@tempcnta}

\def\Vektor(#1,#2)#3(#4,#5){\unskip\leavevmode
  \xcoord#4\relax \ycoord#5\relax
      \raise\ycoord \Einheit\hbox to0pt{\hskip\xcoord \Einheit
         \Vector@(#1,#2){#3}\hss}}

\def\Vector@(#1,#2)#3{\@xarg #1\relax \@yarg #2\relax
\@tempcnta \ifnum\@xarg<0 -\@xarg\else\@xarg\fi
\ifnum\@tempcnta<5\relax
\@linelen=#3\Einheit
\ifnum\@xarg =0 \@vvector
  \else \ifnum\@yarg =0 \@hvector \else \@svector\fi
\fi
\else\@badlinearg\fi}

\def\@hvector{\@hline\hbox to 0pt{\@linefnt
\ifnum \@xarg <0 \@getlarrow(1,0)\hss\else
    \hss\@getrarrow(1,0)\fi}}

\def\@vvector{\ifnum \@yarg <0 \@downvector \else \@upvector \fi}

\def\@svector{\@sline
\@tempcnta\@yarg \ifnum\@tempcnta <0 \@tempcnta=-\@tempcnta\fi
\ifnum\@tempcnta <5
  \hskip -\wd\@linechar
  \@upordown\@clnht \hbox{\@linefnt  \if@negarg
  \@getlarrow(\@xarg,\@yyarg) \else \@getrarrow(\@xarg,\@yyarg) \fi}%
\else\@badlinearg\fi}

\def\@upline{\hbox to \z@{\hskip -.5\Pfadd@cke \vrule width \Pfadd@cke
   height \@linelen depth \z@\hss}}

\def\@downline{\hbox to \z@{\hskip -.5\Pfadd@cke \vrule width \Pfadd@cke
   height \z@ depth \@linelen \hss}}

\def\@upvector{\@upline\setbox\@tempboxa\hbox{\@linefnt\char'66}\raise
     \@linelen \hbox to\z@{\lower \ht\@tempboxa\box\@tempboxa\hss}}

\def\@downvector{\@downline\lower \@linelen
      \hbox to \z@{\@linefnt\char'77\hss}}

\def\@getlarrow(#1,#2){\ifnum #2 =\z@ \@tempcnta='33\else
\@tempcnta=#1\relax\multiply\@tempcnta \sixt@@n \advance\@tempcnta
-9 \@tempcntb=#2\relax\multiply\@tempcntb \tw@
\ifnum \@tempcntb >0 \advance\@tempcnta \@tempcntb\relax
\else\advance\@tempcnta -\@tempcntb\advance\@tempcnta 64
\fi\fi\char\@tempcnta}

\def\@getrarrow(#1,#2){\@tempcntb=#2\relax
\ifnum\@tempcntb < 0 \@tempcntb=-\@tempcntb\relax\fi
\ifcase \@tempcntb\relax \@tempcnta='55 \or
\ifnum #1<3 \@tempcnta=#1\relax\multiply\@tempcnta
24 \advance\@tempcnta -6 \else \ifnum #1=3 \@tempcnta=49
\else\@tempcnta=58 \fi\fi\or
\ifnum #1<3 \@tempcnta=#1\relax\multiply\@tempcnta
24 \advance\@tempcnta -3 \else \@tempcnta=51\fi\or
\@tempcnta=#1\relax\multiply\@tempcnta
\sixt@@n \advance\@tempcnta -\tw@ \else
\@tempcnta=#1\relax\multiply\@tempcnta
\sixt@@n \advance\@tempcnta 7 \fi\ifnum #2<0 \advance\@tempcnta 64 \fi
\char\@tempcnta}

\def\Diagonale(#1,#2)#3{\unskip\leavevmode
  \xcoord#1\relax \ycoord#2\relax
      \raise\ycoord \Einheit\hbox to0pt{\hskip\xcoord \Einheit
         \Line@(1,1){#3}\hss}}
\def\AntiDiagonale(#1,#2)#3{\unskip\leavevmode
  \xcoord#1\relax \ycoord#2\relax 
      \raise\ycoord \Einheit\hbox to0pt{\hskip\xcoord \Einheit
         \Line@(1,-1){#3}\hss}}
\def\Pfad(#1,#2),#3\endPfad{\unskip\leavevmode
  \xcoord#1 \ycoord#2 \thicklines\ZeichnePfad#3\endPfad\thinlines}
\def\ZeichnePfad#1{\ifx#1\endPfad\let\next\relax
  \else\let\next\ZeichnePfad
    \ifnum#1=1
      \raise\ycoord \Einheit\hbox to0pt{\hskip\xcoord \Einheit
         \vrule height\Pfadd@cke width1 \Einheit depth\Pfadd@cke\hss}%
      \advance\xcoord by 1
    \else\ifnum#1=2
      \raise\ycoord \Einheit\hbox to0pt{\hskip\xcoord \Einheit
        \hbox{\hskip-\PfadD@cke\vrule height1 \Einheit width\PfadD@cke depth0pt}\hss}%
      \advance\ycoord by 1
    \else\ifnum#1=3
      \raise\ycoord \Einheit\hbox to0pt{\hskip\xcoord \Einheit
         \Line@(1,1){1}\hss}
      \advance\xcoord by 1
      \advance\ycoord by 1
    \else\ifnum#1=4
      \raise\ycoord \Einheit\hbox to0pt{\hskip\xcoord \Einheit
         \Line@(1,-1){1}\hss}
      \advance\xcoord by 1
      \advance\ycoord by -1
    \fi\fi\fi\fi
  \fi\next}
\def\hSSchritt{\leavevmode\raise-.4pt\hbox to0pt{\hss.\hss}\hskip.2\Einheit
  \raise-.4pt\hbox to0pt{\hss.\hss}\hskip.2\Einheit
  \raise-.4pt\hbox to0pt{\hss.\hss}\hskip.2\Einheit
  \raise-.4pt\hbox to0pt{\hss.\hss}\hskip.2\Einheit
  \raise-.4pt\hbox to0pt{\hss.\hss}\hskip.2\Einheit}
\def\vSSchritt{\vbox{\baselineskip.2\Einheit\lineskiplimit0pt
\hbox{.}\hbox{.}\hbox{.}\hbox{.}\hbox{.}}}
\def\DSSchritt{\leavevmode\raise-.4pt\hbox to0pt{%
  \hbox to0pt{\hss.\hss}\hskip.2\Einheit
  \raise.2\Einheit\hbox to0pt{\hss.\hss}\hskip.2\Einheit
  \raise.4\Einheit\hbox to0pt{\hss.\hss}\hskip.2\Einheit
  \raise.6\Einheit\hbox to0pt{\hss.\hss}\hskip.2\Einheit
  \raise.8\Einheit\hbox to0pt{\hss.\hss}\hss}}
\def\dSSchritt{\leavevmode\raise-.4pt\hbox to0pt{%
  \hbox to0pt{\hss.\hss}\hskip.2\Einheit
  \raise-.2\Einheit\hbox to0pt{\hss.\hss}\hskip.2\Einheit
  \raise-.4\Einheit\hbox to0pt{\hss.\hss}\hskip.2\Einheit
  \raise-.6\Einheit\hbox to0pt{\hss.\hss}\hskip.2\Einheit
  \raise-.8\Einheit\hbox to0pt{\hss.\hss}\hss}}
\def\SPfad(#1,#2),#3\endSPfad{\unskip\leavevmode
  \xcoord#1 \ycoord#2 \ZeichneSPfad#3\endSPfad}
\def\ZeichneSPfad#1{\ifx#1\endSPfad\let\next\relax
  \else\let\next\ZeichneSPfad
    \ifnum#1=1
      \raise\ycoord \Einheit\hbox to0pt{\hskip\xcoord \Einheit
         \hSSchritt\hss}%
      \advance\xcoord by 1
    \else\ifnum#1=2
      \raise\ycoord \Einheit\hbox to0pt{\hskip\xcoord \Einheit
        \hbox{\hskip-2pt \vSSchritt}\hss}%
      \advance\ycoord by 1
    \else\ifnum#1=3
      \raise\ycoord \Einheit\hbox to0pt{\hskip\xcoord \Einheit
         \DSSchritt\hss}
      \advance\xcoord by 1
      \advance\ycoord by 1
    \else\ifnum#1=4
      \raise\ycoord \Einheit\hbox to0pt{\hskip\xcoord \Einheit
         \dSSchritt\hss}
      \advance\xcoord by 1
      \advance\ycoord by -1
    \fi\fi\fi\fi
  \fi\next}
\def\Koordinatenachsen(#1,#2){\unskip
 \hbox to0pt{\hskip-.5pt\vrule height#2 \Einheit width.5pt depth1 \Einheit}%
 \hbox to0pt{\hskip-1 \Einheit \xcoord#1 \advance\xcoord by1
    \vrule height0.25pt width\xcoord \Einheit depth0.25pt\hss}}
\def\Koordinatenachsen(#1,#2)(#3,#4){\unskip
 \hbox to0pt{\hskip-.5pt \ycoord-#4 \advance\ycoord by1
    \vrule height#2 \Einheit width.5pt depth\ycoord \Einheit}%
 \hbox to0pt{\hskip-1 \Einheit \hskip#3\Einheit
    \xcoord#1 \advance\xcoord by1 \advance\xcoord by-#3
    \vrule height0.25pt width\xcoord \Einheit depth0.25pt\hss}}
\def\Gitter(#1,#2){\unskip \xcoord0 \ycoord0 \leavevmode
  \LOOP\ifnum\ycoord<#2
    \loop\ifnum\xcoord<#1
      \raise\ycoord \Einheit\hbox to0pt{\hskip\xcoord \Einheit\Punkt\hss}%
      \advance\xcoord by1
    \repeat
    \xcoord0
    \advance\ycoord by1
  \REPEAT}
\def\Gitter(#1,#2)(#3,#4){\unskip \xcoord#3 \ycoord#4 \leavevmode
  \LOOP\ifnum\ycoord<#2
    \loop\ifnum\xcoord<#1
      \raise\ycoord \Einheit\hbox to0pt{\hskip\xcoord \Einheit\Punkt\hss}%
      \advance\xcoord by1
    \repeat
    \xcoord#3
    \advance\ycoord by1
  \REPEAT}
\def\Label#1#2(#3,#4){\unskip \xdim#3 \Einheit \ydim#4 \Einheit
  \def\lo{\advance\xdim by-.5 \Einheit \advance\ydim by.5 \Einheit}%
  \def\llo{\advance\xdim by-.25cm \advance\ydim by.5 \Einheit}%
  \def\loo{\advance\xdim by-.5 \Einheit \advance\ydim by.25cm}%
  \def\o{\advance\ydim by.25cm}%
  \def\ro{\advance\xdim by.5 \Einheit \advance\ydim by.5 \Einheit}%
  \def\rro{\advance\xdim by.25cm \advance\ydim by.5 \Einheit}%
  \def\roo{\advance\xdim by.5 \Einheit \advance\ydim by.25cm}%
  \def\l{\advance\xdim by-.30cm}%
  \def\r{\advance\xdim by.30cm}%
  \def\lu{\advance\xdim by-.5 \Einheit \advance\ydim by-.6 \Einheit}%
  \def\llu{\advance\xdim by-.25cm \advance\ydim by-.6 \Einheit}%
  \def\luu{\advance\xdim by-.5 \Einheit \advance\ydim by-.30cm}%
  \def\u{\advance\ydim by-.30cm}%
  \def\ru{\advance\xdim by.5 \Einheit \advance\ydim by-.6 \Einheit}%
  \def\rru{\advance\xdim by.25cm \advance\ydim by-.6 \Einheit}%
  \def\ruu{\advance\xdim by.5 \Einheit \advance\ydim by-.30cm}%
  #1\raise\ydim\hbox to0pt{\hskip\xdim
     \vbox to0pt{\vss\hbox to0pt{\hss$#2$\hss}\vss}\hss}%
}

\font@\twelverm=cmr10 scaled\magstep1
\font@\twelveit=cmti10 scaled\magstep1
\font@\twelvebf=cmbx10 scaled\magstep1
\font@\twelvei=cmmi10 scaled\magstep1
\font@\twelvesy=cmsy10 scaled\magstep1
\font@\twelveex=cmex10 scaled\magstep1

\newtoks\twelvepoint@
\def\twelvepoint{\normalbaselineskip15\p@
 \abovedisplayskip15\p@ plus3.6\p@ minus10.8\p@
 \belowdisplayskip\abovedisplayskip
 \abovedisplayshortskip\z@ plus3.6\p@
 \belowdisplayshortskip8.4\p@ plus3.6\p@ minus4.8\p@
 \textonlyfont@\rm\twelverm \textonlyfont@\it\twelveit
 \textonlyfont@\sl\twelvesl \textonlyfont@\bf\twelvebf
 \textonlyfont@\smc\twelvesmc \textonlyfont@\tt\twelvett
%
 \ifsyntax@ \def\big##1{{\hbox{$\left##1\right.$}}}%
  \let\Big\big \let\bigg\big \let\Bigg\big
 \else
  \textfont\z@=\twelverm  \scriptfont\z@=\tenrm  \scriptscriptfont\z@=\sevenrm
  \textfont\@ne=\twelvei  \scriptfont\@ne=\teni  \scriptscriptfont\@ne=\seveni
  \textfont\tw@=\twelvesy \scriptfont\tw@=\tensy \scriptscriptfont\tw@=\sevensy
  \textfont\thr@@=\twelveex \scriptfont\thr@@=\tenex
        \scriptscriptfont\thr@@=\tenex
  \textfont\itfam=\twelveit \scriptfont\itfam=\tenit
        \scriptscriptfont\itfam=\tenit
  \textfont\bffam=\twelvebf \scriptfont\bffam=\tenbf
        \scriptscriptfont\bffam=\sevenbf
  \setbox\strutbox\hbox{\vrule height10.2\p@ depth4.2\p@ width\z@}%
  \setbox\strutbox@\hbox{\lower.6\normallineskiplimit\vbox{%
        \kern-\normallineskiplimit\copy\strutbox}}%
 \setbox\z@\vbox{\hbox{$($}\kern\z@}\bigsize@=1.4\ht\z@
 \fi
 \normalbaselines\rm\ex@.2326ex\jot3.6\ex@\the\twelvepoint@}

\font@\fourteenrm=cmr10 scaled\magstep2
\font@\fourteenit=cmti10 scaled\magstep2
\font@\fourteensl=cmsl10 scaled\magstep2
\font@\fourteensmc=cmcsc10 scaled\magstep2
\font@\fourteentt=cmtt10 scaled\magstep2
\font@\fourteenbf=cmbx10 scaled\magstep2
\font@\fourteeni=cmmi10 scaled\magstep2
\font@\fourteensy=cmsy10 scaled\magstep2
\font@\fourteenex=cmex10 scaled\magstep2
\font@\fourteenmsa=msam10 scaled\magstep2
\font@\fourteeneufm=eufm10 scaled\magstep2
\font@\fourteenmsb=msbm10 scaled\magstep2
\newtoks\fourteenpoint@
\def\fourteenpoint{\normalbaselineskip15\p@
 \abovedisplayskip18\p@ plus4.3\p@ minus12.9\p@
 \belowdisplayskip\abovedisplayskip
 \abovedisplayshortskip\z@ plus4.3\p@
 \belowdisplayshortskip10.1\p@ plus4.3\p@ minus5.8\p@
 \textonlyfont@\rm\fourteenrm \textonlyfont@\it\fourteenit
 \textonlyfont@\sl\fourteensl \textonlyfont@\bf\fourteenbf
 \textonlyfont@\smc\fourteensmc \textonlyfont@\tt\fourteentt
%
 \ifsyntax@ \def\big##1{{\hbox{$\left##1\right.$}}}%
  \let\Big\big \let\bigg\big \let\Bigg\big
 \else
  \textfont\z@=\fourteenrm  \scriptfont\z@=\twelverm  \scriptscriptfont\z@=\tenrm
  \textfont\@ne=\fourteeni  \scriptfont\@ne=\twelvei  \scriptscriptfont\@ne=\teni
  \textfont\tw@=\fourteensy \scriptfont\tw@=\twelvesy \scriptscriptfont\tw@=\tensy
  \textfont\thr@@=\fourteenex \scriptfont\thr@@=\twelveex
        \scriptscriptfont\thr@@=\twelveex
  \textfont\itfam=\fourteenit \scriptfont\itfam=\twelveit
        \scriptscriptfont\itfam=\twelveit
  \textfont\bffam=\fourteenbf \scriptfont\bffam=\twelvebf
        \scriptscriptfont\bffam=\tenbf
  \setbox\strutbox\hbox{\vrule height12.2\p@ depth5\p@ width\z@}%
  \setbox\strutbox@\hbox{\lower.72\normallineskiplimit\vbox{%
        \kern-\normallineskiplimit\copy\strutbox}}%
 \setbox\z@\vbox{\hbox{$($}\kern\z@}\bigsize@=1.7\ht\z@
 \fi
 \normalbaselines\rm\ex@.2326ex\jot4.3\ex@\the\fourteenpoint@}

\catcode`\@=13
\def\[{\left[}
\def\]{\right]}
\define\twoline#1#2{\line{\hfill{\smc #1}\hfill{\smc #2}\hfill}}

\def\mypic#1{\epsffile{#1}}

\topmatter
\title The number of centered lozenge tilings of a symmetric hexagon
\endtitle
\author M.~Ciucu%
\footnote"$^\dagger$"{\hbox{Research supported in part by the
MSRI, Berkeley.}} and C.~Krattenthaler%
$^\dagger$
\endauthor
\affil School of Mathematics,\\
Institute for Advanced Study, \\
Princeton, NJ 08540, USA\\
email: ciucu\@math.ias.edu\\\vskip6pt
Institut f\"ur Mathematik der Universit\"at Wien,\\
Strudlhofgasse 4, A-1090 Wien, Austria.\\
e-mail: KRATT\@Pap.Univie.Ac.At\\
WWW: \tt http://radon.mat.univie.ac.at/People/kratt
\endaffil
\address School of Mathematics,
Institute for Advanced Study,
Princeton, NJ 08540, USA.
\endaddress
\address Institut f\"ur Mathematik der Universit\"at Wien,
Strudlhofgasse 4, A-1090 Wien, Austria.
\endaddress
\abstract
Propp conjectured \cite{\PropAA} that the number of
lozenge tilings of a semiregular hexagon of sides $2n-1$, $2n-1$ and $2n$
which contain the central unit rhombus is precisely one third of the total
number of lozenge tilings. Motivated by this, we consider the more general
situation of a semiregular hexagon of sides $a$, $a$ and $b$. We prove
explicit formulas for the number of lozenge tilings of these hexagons
containing the central unit rhombus, and obtain Propp's conjecture as a
corollary of our results. 
\endabstract
\endtopmatter
\document

\leftheadtext{M. Ciucu and C. Krattenthaler}
\rightheadtext{The number of centered lozenge tilings}

\subhead 1. Introduction\endsubhead
Let $a$, $b$ and $c$ be positive integers, and consider a semiregular
hexagon of sides $a$, $b$ and $c$ (i.e., all angles have 120 degrees
and the sides have, in order, lengths $a$, $b$, $c$, $a$, $b$, $c$).
By a well-known bijection \cite{\DT}, the number of tilings of this hexagon
by rhombi of unit edge-length and angles of 60 and 120 degrees
(we call such a rhombus a 
{\it lozenge} and such tilings {\it lozenge tilings}) is equal
to the number $P(a,b,c)$ of plane
partitions contained in an $a\times b\times c$ box. In turn,
by a famous result of MacMahon \cite{\MM}, the latter is given by the product
  
$$P(a,b,c)=\prod _{i=1} ^{a}\prod _{j=1} ^{b}\prod _{k=1} ^{c}\frac {i+j+k-1}
{i+j+k-2}.\tag\Aa$$

The starting point of this paper is a conjecture of Propp \cite{\PropAA}
stating that for a semiregular hexagon of sides $2n-1$, $2n-1$ and $2n$,
precisely one third of its lozenge tilings contain the central lozenge.
Call the lozenge tilings with this property {\it centered}.
In this paper we consider the following more general problem: for a
semiregular hexagon of sides $a$, $a$ and $b$, how many of its tilings
are centered? It is easy to see that such a hexagon has a
central lozenge only if $a$ and $b$ have opposite parity. The two cases are
addressed in Theorems 1 and 2 below. 

For any nonnegative integer $m$ and positive integer $n$ define

$$Q(m,n)=\frac {(2n)!^2\,(2m)!\,(m+2n-1)!} {2\cdot n!^2\,m!\,(2m+4n-2)!}
\bigg( \sum _{i=0} ^{n-1}\frac {(-1)^{n-i-1}} {(2n-2i-1)}\frac
{(m+n-i)_{2i}} {i!^2}\bigg),\tag\AC$$
where the shifted factorial $(a)_k$ is defined by $(a)_k:=a(a+1)\cdots(a+k-1)$,
$k\ge1$, and $(a)_0:=1$.

\proclaim{Theorem~\TA}Let $m$ be a nonnegative integer and $n$ a positive integer.
The number of centered lozenge tilings of a semiregular
hexagon with sides $2n-1$, $2n-1$ and $2m$ is $Q(m,n)P(2n-1,2n-1,2m)$.
\endproclaim

\proclaim{Theorem~\TB}Let $m$ and $n$ be positive integers.
The number of centered lozenge tilings of a semiregular
hexagon with sides $2n$, $2n$ and $2m-1$ is $Q(m,n)P(2n,2n,2m-1)$.
\endproclaim

In the case when $m$ equals $n$, the expression $Q(m,n)$ evaluates to 1/3
(this is due to a remarkable simplification of the sum in (\AC) in this case).
Thus, the statement in Propp's conjecture follows from Theorem 1.

\proclaim{Corollary~\TC}Let $n$ be a positive integer.
Exactly one third of the lozenge tilings of a
semiregular hexagon with sides $2n-1$, $2n-1$ and $2n$ contain the
central lozenge. The same is true for a semiregular hexagon with
sides $2n$, $2n$ and $2n-1$.
\endproclaim

Theorems~\TA\ and \TB\ have also been found (in an equivalent form)
by Helfgott and Gessel \cite{\HeGeAA, Theorem~2}, using a completely
different method. They also obtained Corollary~\TC\ as a corollary of
their results.

\medskip
For $m\neq n$ the sum in (\AC) does not seem to simplify. However, if $m$
and $n$ approach infinity so that their ratio approaches some non-negative
real number $a$, $Q(m,n)$ turns out to
approach the value $\frac {2} {\pi}\arcsin(1/(a+1))$.

\proclaim{Corollary~\TD}Let $a$ be any nonnegative real number.
For $m\sim an$, the proportion
of the rhombus tilings that contain the
central rhombus in the total number of rhombus tilings of a
semiregular hexagon with sides $2n-1$, $2n-1$ and $2m$ is
$\sim\ \frac {2} {\pi}\arcsin(1/(a+1))$
as $n$ tends to infinity.
The same is true for a
semiregular hexagon with sides $2n$, $2n$ and $2m-1$.
\endproclaim

\flushpar
{\it Remark.} Using the bijection \cite{\DT} between lozenge tilings and
plane partitions, the statement of Theorem 1 can be interpreted as follows.
Let $\Cal P$ be the set of plane partitions $(a_{ij})$ of square shape
$(2n-1)^{2n-1}$, with entries between 0 and $2m$. Let ${\Cal P}_k$ be the
subset consisting of the plane partitions for which $a_{n+k,n+k}=m+k$, for
$-\min(n-1,m)\leq k\leq\min(n-1,m)$. Then the number of elements in
the union of the ${\Cal P}_k$'s is $Q(m,n)P(2n-1,2n-1,2m)$ (this union is clearly
disjoint).   

A similar interpretation can be given to the statement of Theorem 2.

\smallskip
The rest of the paper is devoted to giving proofs of Theorems~\TA\ and
\TB, and of Corollaries~\TC\ and \TD. In Section~2 we provide proofs
of Corollaries~\TC\ and \TD, and outline the proofs of
Theorems~\TA\ and \TB, the latter consisting of several steps.
The details of these steps are then given in detail in the subsequent
sections. These steps are the following. First,
an application of the first author's Matchings Factorization Theorem
\cite{\CiucAB, Theorem~1.2} allows to reduce our problem to the
enumeration of lozenge tilings 
of {\it simply}-connected regions. This is described in
Section~3. Then, in Section~4, lozenge tilings are translated into
nonintersecting lattice paths. By the main theorem on nonintersecting
lattice paths \cite{\GeViAA, \GeViAB}, the number(s) of nonintersecting
lattice paths that we are interested in
can be immediately written down in form of a
determinant (see Lemmas~\TK, \TKa, \TL). Finally, in Section~5, these
determinants are evaluated (see Lemmas~\TN\ and \TO).

\subhead 2. Outline of proofs\endsubhead
Here we outline the proofs of Theorems~\TA\ and \TB\ and we deduce
Corollaries~\TC\ and \TD. We fill in the details in the subsequent sections.
Denote by $L(R)$ the number of lozenge tilings of the region $R$.

\smallskip
{\smc Proof of Theorem~\TA}.
In Section~3 it is shown that
the number of centered lozenge tilings of a semiregular
hexagon with sides $2n-1$, $2n-1$ and $2m$ equals $2^{2n-2}$ times the product
of the number of lozenge tilings of two
regions of the triangular lattice, $H^+$ and $H^-$ (see (3.2) and Figure~3.3).
Then, in Section~4 we use the Gessel-Viennot method of 
nonintersecting lattice paths to obtain determinantal expressions for
$L(H^+)$ (see Lemma~\TK) and $L(H^-)$ (see Lemma~\TL).
Finally, in Section~5 we evaluate these determinants (see Lemmas~\TN\
(with $N=2n-2$) and \TO). After some manipulation of the expressions on the 
right-hand sides of (5.2) and (5.3) one obtains the statement of the Theorem.
\quad \quad \qed

\smallskip
{\smc Proof of Theorem~\TB}. We proceed along the same lines as in the proof of
Theorem 1. In Section~3 we show that the number of centered lozenge
tilings of a
semiregular hexagon with sides $2n$, $2n$ and $2m-1$ equals $2^{2n-1}$ times
the product of the number of lozenge tilings of two
regions of the triangular lattice, $\bar H^+$ and $H^-$ (see (3.3)
and Figure~3.4). The
region $H^-$ is the same as the one appearing in the proof of Theorem 1, while
$\bar H^+$ differs from $H^+$ only in the sizes of its sides. The
determinant evaluations in Lemmas~\TN\ (with $N=2n$ and $m$ replaced by $m-1$)
and 16 lead, after some manipulation of the expressions involved, to the statement
of the Theorem. \quad \quad \qed

\smallskip
{\smc Proof of Corollary~\TC}.
We have to compute the value of the expression on the right-hand side of (\AC) for
$m=n$. Clearly, except for
trivial manipulations, we will be done once we are able to evaluate
the sum in (\AC) for $m=n$.

We claim that
$$\sum _{i=0} ^{n-1}\frac {(-1)^{n-i-1}} {(2n-2i-1)}\frac
{(2n-i)_{2i}} {i!^2}=3^{n-1}\frac {\prod _{i=1} ^{n-1}(6i-1)(6i+1)}
{(2n-1)!!^2}\tag\BA$$
(where the empty product is defined to be 1).
Let us denote the sum by $S(n)$ and its summand by $F(n,i)$.
We use the Gosper--Zeilberger algorithm \cite{\PeWZAA, \ZeilAM, \ZeilAV}
to obtain the relation
$$   n {{\left( 2 n +1\right) }^2} F(n+1,i) -
 3n\left( 6 n-1 \right)   \left(  6 n+1 \right)  F(n,i)
=G(n,i+1)-G(n,i),\tag\BB$$
with
$$\multline
G(n,i)={{\left( -1 \right) }^{n-i}}
     \left( -3 + 9\, i - 6\, {i^2} - 30 \,n + 62 \,i n - 28 \,{i^2} n -
       104 \,{n^2} + 104 \,i {n^2} - 112 \,{n^3} \right)\\
\times{{  i^2\,
     ({ \textstyle 2 n- i + 2 }) _{2i-2} }\over
   {\left( 2n - 2 i + 1 \right) \, {{i  !}^2}}}.
\endmultline$$
Summation of the relation (\BB) from $i=0$ to $i=n$, little
rearrangement, and division by $n$ on both sides,
leads to the recurrence
$$    {{\left( 2 n+1 \right) }^2} S(n+1) -
 3\left( 6 n-1 \right)  \left(  6 n+1 \right)  S(n)=0
$$
for the sum in (\BA). (Paule and Schorn's \cite{\PaScAA} {\sl Mathematica}
implementation of the Gosper--Zeilberger algorithm,
which is the one we used, gives this recurrence directly.) Since
$S(1)=1$, and since the right-hand side of (\BA) satisfies the same
recurrence, equation (\BA) is proved, and, thus, the Corollary also.
\quad \quad \qed

{\smc Proof of Corollary~\TD}.
First let $a>0$. We have to determine the limit of $Q(m,n)$ as $n$ tends
to infinity, and where the relation
between $m$ and $n$ is fixed by $m\sim an$. Clearly, the ``difficult"
part of this asymptotic computation is to find the asymptotics of the
sum in (\AC). It turns out that it is convenient to manipulate this
sum first, before taking the limit $n\to\infty$. We reverse the order
of summation in the sum in (\AC), and then are able to rewrite the sum
using the standard hypergeometric notation
$${}_r F_s\!\left[\matrix a_1,\dots,a_r\\ b_1,\dots,b_s\endmatrix;
z\right]=\sum _{k=0} ^{\infty}\frac {\po{a_1}{k}\cdots\po{a_r}{k}}
{k!\,\po{b_1}{k}\cdots\po{b_s}{k}} z^k\ .\tag\Ac$$
(The reader should notice that the sum in
(\AC) cannot be directly converted into a hypergeometric series since
the upper bound on the summation index, $n-1$, is ``artificial",
i.e., the sum is changed if we extend the range of summation to all
nonnegative $i$. If however the order of summation is reversed, this
problem is removed because of the presence of the term $i!$ in the
denominator of the summand.) Thus, expression (\AC) is converted into
$$
\frac {(2n)!^2\,(2m)!\,(m+2n-1)!} {2\cdot n!^2\,m!\,(2m+4n-2)!}
{{     ({ \textstyle m+1}) _{2n-2} }\over {{{\left( n-1 \right) !}^2}}
   }
{} _{4} F _{3} \!\left [ \matrix { 1, {1\over 2}, 1 - n, 1 - n}\\ { 1 + m, 2
      - m - 2\,n, {3\over 2}}\endmatrix ; {\displaystyle 1}\right ].
$$
Next we apply Bailey's transformation for a balanced $_4F_3$-series
(see \cite{\SlatAC, (4.3.5.1)}),
$$\multline
{} _{4} F _{3} \!\left [ \matrix { a, b, c, -N}\\ { e, f, 1 + a + b + c - e -
   f - N}\endmatrix ; {\displaystyle 1}\right ]  \\{{( e-a)_N\, (f-a) _{N}}\over  {(e)_N\,( f) _{N}}}
  {} _{4} F _{3} \!\left [ \matrix { -N, a, 1 + a + c - e - f - N, 1 + a + b -
    e - f - N}\\ { 1 + a + b + c - e - f - N, 1 + a - e - N, 1 + a - f -
    N}\endmatrix ; {\displaystyle 1}\right ]   ,
\endmultline$$
where $N$ is a nonnegative integer. Thus we obtain the expression
$$
\frac {(2n)!^2\,(2m)!\,(m+2n-1)!} {2\cdot n!^2\,m!\,(2m+4n-2)!}
{{    ({ \textstyle m}) _{n-1} \,({ \textstyle m + n+1}) _{n-1} }\over
   {{{\left( n-1 \right) !}^2}}}
{} _{4} F _{3} \!\left [ \matrix { 1 - n, 1, 1, {1\over 2} + n}\\ { {3\over
      2}, 2 - m - n, 1 + m + n}\endmatrix ; {\displaystyle 1}\right ].
$$
Now we substitute $m\sim an$ and perform the limit $n\to\infty$.
Using Stirling's formula, it is easy to determine the limit for the
quotient in front of the $_4F_3$-series. It is
$$\frac {2\sqrt{a(a+2)}} {\pi (a+1)^2}.\tag\AE$$
For the $_4F_3$-series itself, we may exchange limit and summation,
because of uniform convergence. This gives
$$\lim_{n\to\infty}
{} _{4} F _{3} \!\left [ \matrix { 1 - n, 1, 1, {1\over 2} + n}\\ { {3\over
      2}, 2 - m - n, 1 + m + n}\endmatrix ; {\displaystyle 1}\right
]={} _{2} F _{1} \!\left [ \matrix { 1, 1}\\ { {3\over
      2}}\endmatrix ; {\displaystyle \frac {1} {(a+1)^2}}\right ].
\tag\AF$$
(Recall that $m\sim an$.)
Combining (\AE) and (\AF), and using the
identity (see \cite{\PrBMAA, p. 463, (133)})
$${}_2F_1\!\[\matrix 1,1\\\frac {3} {2}\endmatrix; z\]=\frac {\arcsin\sqrt{z}}
{\sqrt{z(1-z)}}$$
in (\AF), we obtain the desired limit $\frac {2}
{\pi}\arcsin(1/(a+1))$.

Finally we address the case $a=0$, which means that $m\sim 0$. Then
the proportion (\AC) is arbitrarily close to the expression which
results from (\AC) when $m\to0$. In particular,
according to Theorem~\TA, this expression gives the proportion
of the rhombus tilings that contain the
central rhombus in the total number of rhombus tilings of the
semiregular hexagon with sides $2n-1$, $2n-1$ and $0$. That
proportion is simply 1 since, trivially, there is exactly one such
rhombus tiling, and it does contain the central rhombus. The value of
1 agrees with the claimed expression $\frac {2} {\pi}\arcsin(1/(a+1))$
evaluated at $a=0$.
\quad \quad \qed

\subhead 3. Reduction to simply-connected regions\endsubhead
One useful way to approach certain tiling enumeration problems is to biject
them with nonintersecting lattice paths, and then use the Gessel-Viennot
determinant theorem 
\cite{\GeViAA, \GeViAB}. This approach seems to be especially
appropriate if the entries of the Gessel-Viennot matrix have a simple
expression. In the case of the $(2n-1)\times(2n-1)\times2m$
semiregular hexagon with the
central lozenge removed (whose tilings can clearly be identified with the centered
tilings we are concerned with) this is not quite the case. However, one can get
around this using the Factorization Theorem for perfect matchings presented
in \cite{\CiucAB, Theorem~1.2}.

Consider the tiling of the plane by unit equilateral triangles, illustrated
in Figure 3.1. Define a {\it region} to be the union of finitely many such
unit triangles. Suppose the region $R$ is symmetric with respect to the
horizontal symmetry axis $l$. Suppose further that the unit triangles of $R$
crossed by $l$  can be grouped in pairs such that the two triangles in a
pair share an edge, forming a rhombic tile. Let $T_1,\dotsc,T_k$ be these
rhombi. Let $P$ be the zig-zag lattice path that borders the tiles $T_i$ on
their upper boundary (see Figure 3.2). Define $R^+$ and $R^-$ to be the pieces of $R$
above and below $P$, respectively. Then the Factorization Theorem
of \cite{\CiucAB} implies

\topinsert
\twoline{\mypic{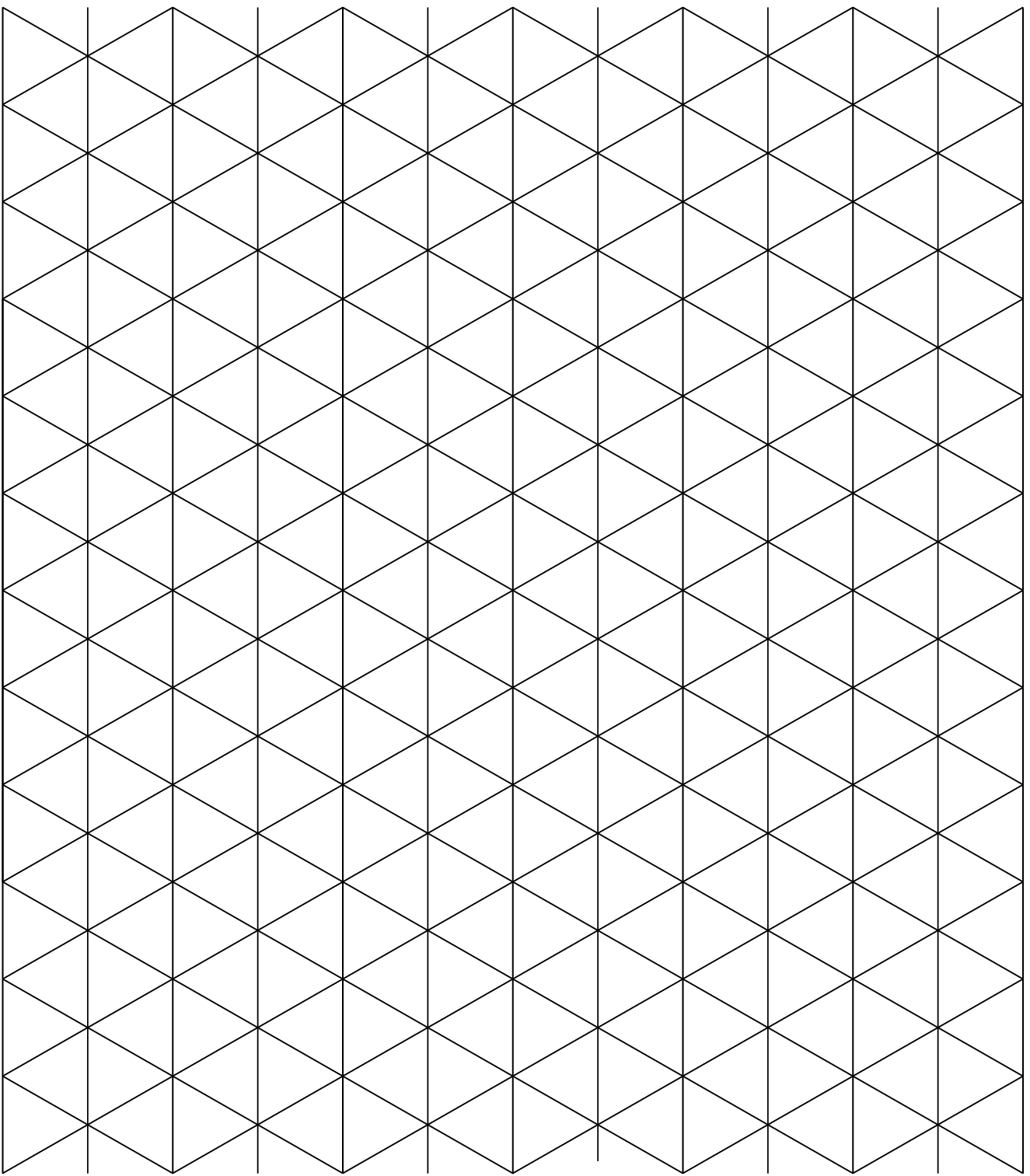}}{\mypic{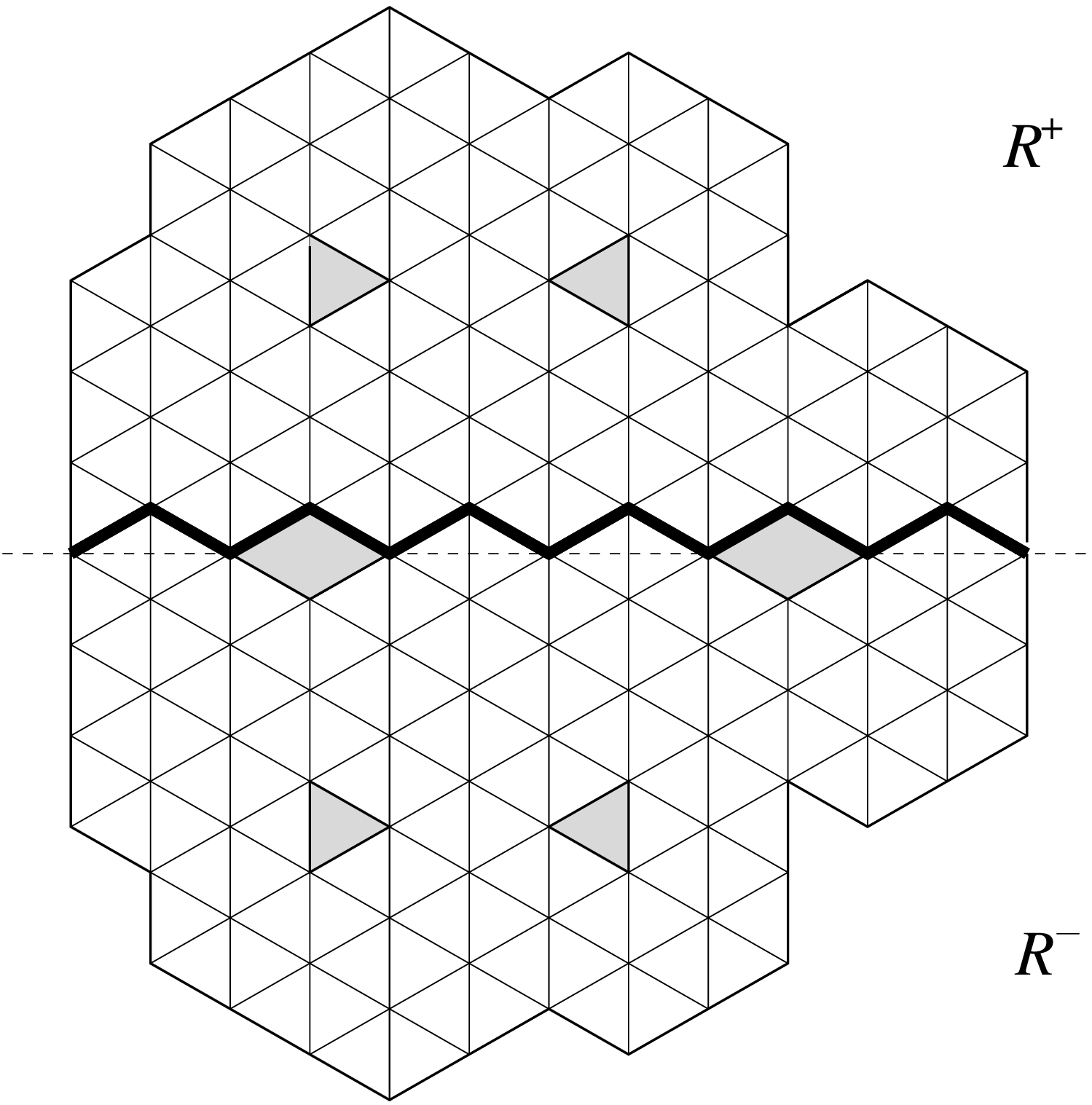}}
\twoline{Figure~3.1}{Figure~3.2}
\endinsert

$$L(R)=2^kL(R^+)L^*(R^-),\tag3.1$$
where $L(R)$ is the number of lozenge tilings of $R$, and $L^*(R^-)$ is the
weighted count of the lozenge tilings of $R^-$ assigning weight $2^{-i}$ to
a lozenge tiling containing $i$ of the rhombi $T_1,\dotsc,T_k$ (a similar
corollary of the Factorization Theorem is given in Remark 2.3 of
\cite{\CiucAB}
for the square lattice).

\topinsert
\twoline{\mypic{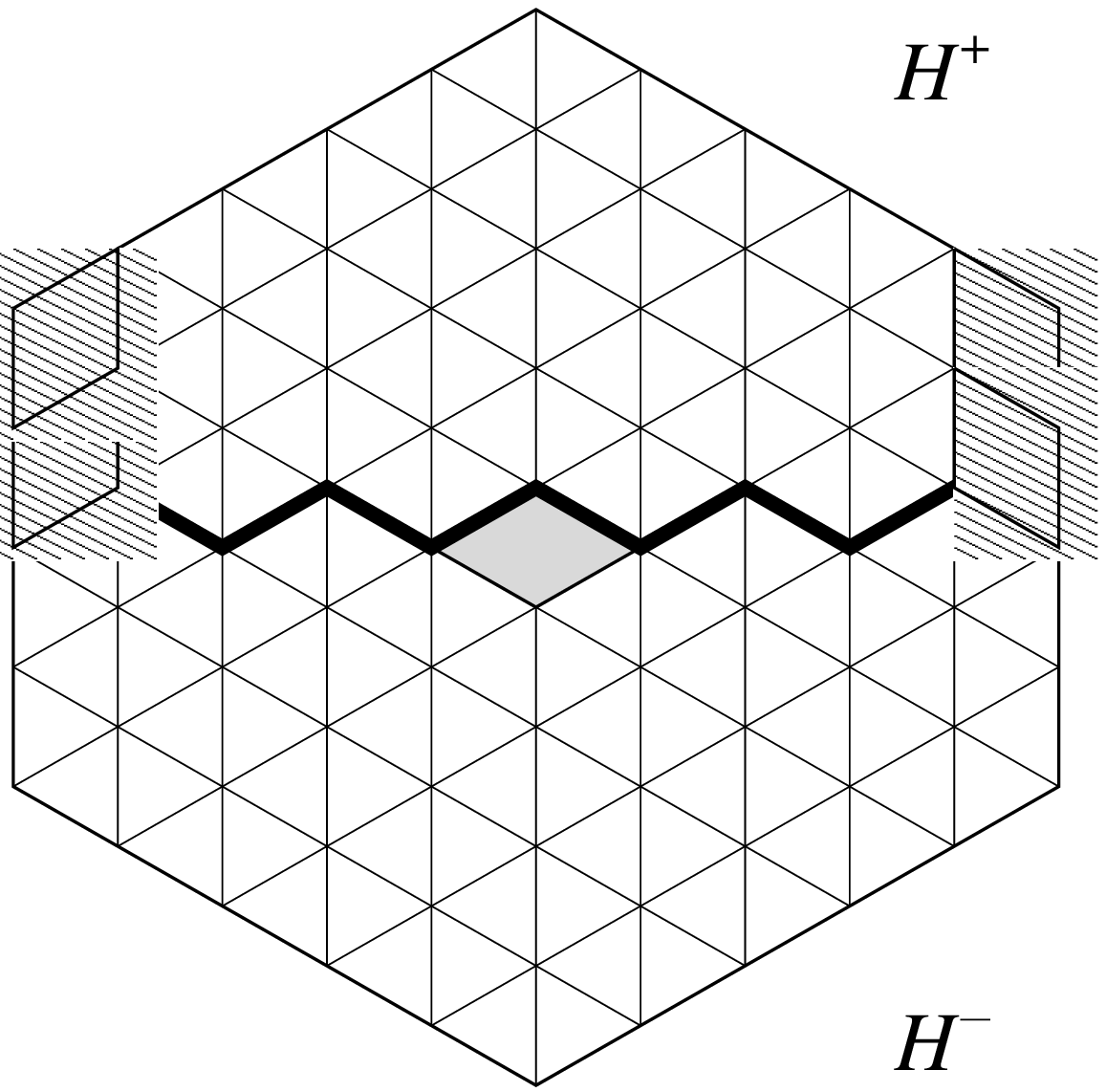}}{\mypic{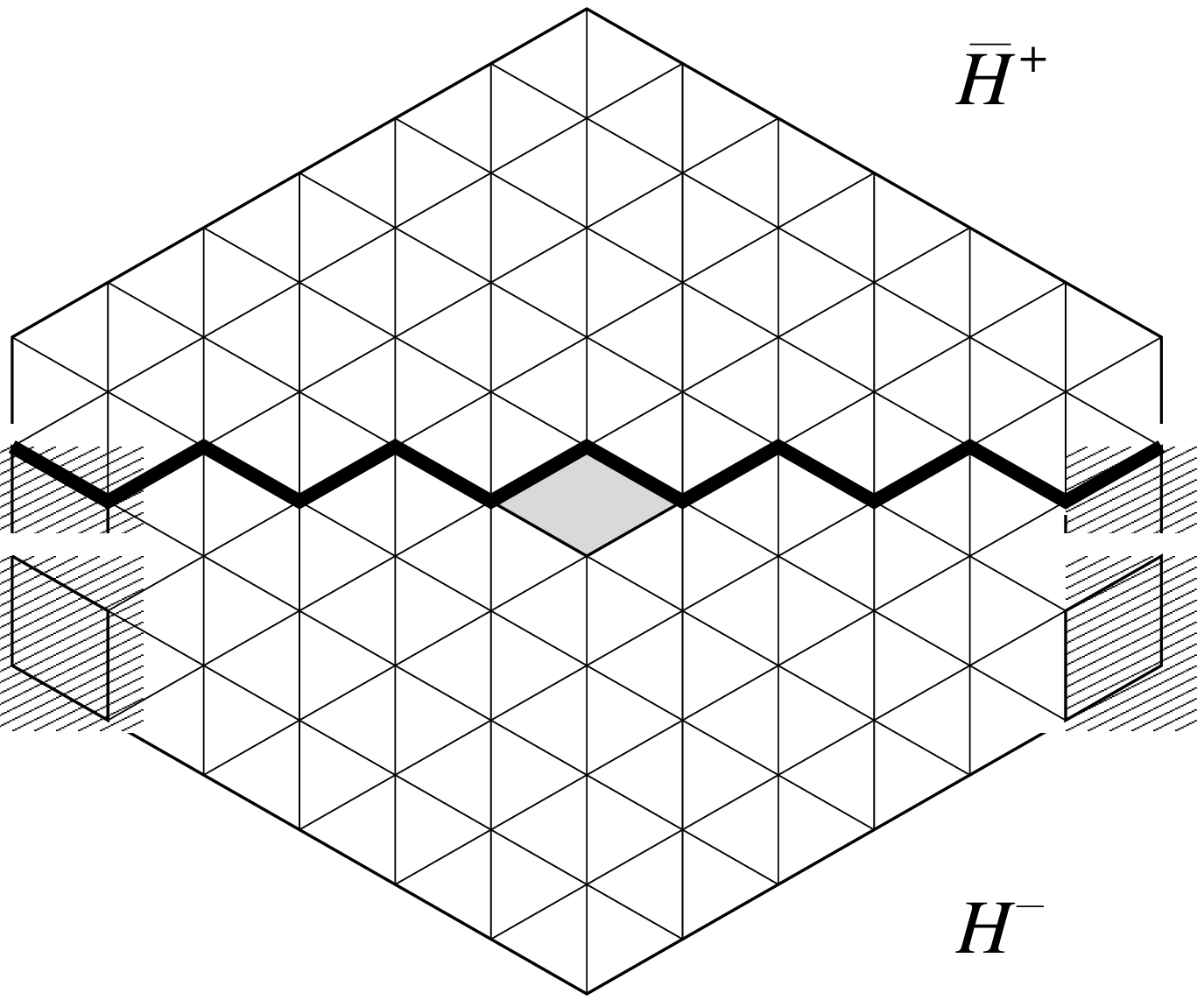}}
\twoline{Figure~3.3}{Figure~3.4}
\endinsert

Let $H$ be the region obtained from a semiregular hexagon with side-lengths
$2n-1,2n-1,2m$ by removing the central lozenge (see Figure 3.3). By (3.1) we obtain
$$L(H)=2^{2n-2}L(H^+)L^*(H^-),\tag3.2$$
where the regions $H^+$ and $H^-$ are indicated in Figure 3.3 ($H^+$ is
obtained from the piece above the zig-zag path by removing $2m$
forced lozenges, which are indicated by a shading). 
Similarly, applying (3.1) to the hexagonal region
$\bar{H}$ (illustrated in Figure 3.4) with side-lengths $2n,2n,2m-1$ and
central lozenge removed we obtain
$$L(\bar{H})=2^{2n-1}L(\bar{H}^+)L^*(H^-).\tag3.3$$
(Indeed, the region obtained from the bottom piece in Figure 3.4 after
removing the forced tiles, which are again indicated by a shading, 
is the same as the region $H^-$ in (3.2); this explains the last factor in (3.3);
$\bar{H}^+$ is shown in Figure 3.4).

\subhead 4. Enumeration of lozenge tilings and nonintersecting
lattice paths\endsubhead
In this section we make use of the standard
encoding of lozenge tilings in terms of non-intersecting lattice paths (Figure 4.1.
illustrates this in the case of the region $H^+$). Thus we transform the problem of 
enumerating lozenge tilings of the regions that arose in Section~2
into the problem of enumerating certain families of
nonintersecting lattice paths. This allows us to derive determinantal formulas
for the number of lozenge tilings we are interested in. 

\topinsert
\centerline{\mypic{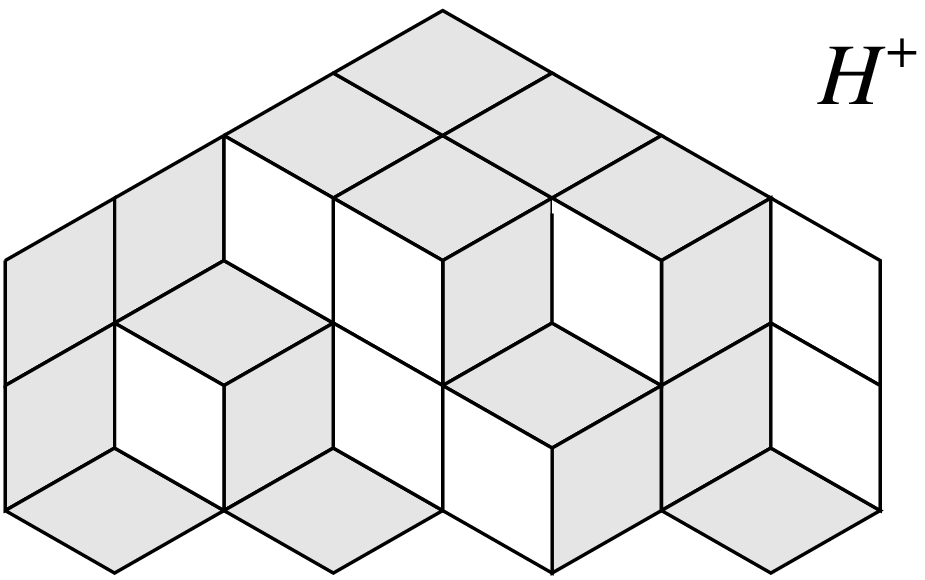}}
\centerline{\smc Figure~4.1}
\endinsert

Our next Lemma exhibits an expression for $L(H^+)$.
\proclaim{Lemma~\TK}We have
$$L(H^+)=\det_{1\le i,j\le 2n-2}\(\binom{2n+m-i-1}{m+i-j}\).\tag\DA$$
\endproclaim
\demo{Proof}
Using the correspondence with nonintersecting lattice paths, one obtains that 
the number of lozenge tilings of the upper ``half hexagon" $H^+$ (see Figure
3.3) equals the
number of families $(P_1,P_2,\dots,P_{2n-2})$
of nonintersecting lattice paths (consisting of horizontal positive unit
steps and vertical negative unit steps)
in which $P_i$ runs from $A_i=(2i,i+m)$ to $E_i=(2n+i-1,i)$,
$i=1,2,\dots,2n-2$. By the main theorem of nonintersecting lattice paths
\cite{\GeViAB, Cor.~2; \StemAE, Theorem~1.2}, this number is given by
the determinant
$$\det_{1\le i,j\le 2n-2}\big(\vert\Cal P(A_i\to E_j)\vert\big),$$
where $\Cal P(A\to E)$ denotes the set of paths from $A$ to $E$ and
$\vert \Cal M\vert$ is the cardinality of the set $\Cal M$.
Obviously, the number $\vert\Cal P(A_i\to E_j)\vert$ of paths from
$A_i$ to $E_j$ equals the binomial $\binom{2n+m-i-1}{m+i-j}$. This
establishes the Lemma.\quad \quad \qed

\enddemo

In the same way, we may derive a determinantal formula for
$L(\bar H^+)$. As we have already noted in Section~2, the only difference
between the regions $H^+$ and $\bar H^+$ is in their side lengths.
To be precise, the vertical sides of $H^+$ have length $m$, while
those of $\bar H^+$ have length $m-1$. On the other hand,
the slanted sides of $H^+$ have length $2n-2$, while
those of $\bar H^+$ have length $2n$. Hence, we will get a
formula for $L(\bar H^+)$ by replacing $m$ by $m-1$ and $n$ by
$n+1$ in the formula for $L(H^+)$.
\proclaim{Lemma~\TKa}We have

\vbox{
$$L(\bar H^+)=\det_{1\le i,j\le
2n}\(\binom{2n+m-i}{m+i-j-1}\).\tag\DAa$$
\line{\hfil\hbox{\qed}\quad \quad }
}
\endproclaim

Next we turn to the number $L(H^-)$. In analogy to the preceding, we
derive also a determinantal expression for $L(H^-)$, using
nonintersecting lattice paths. However, the resulting determinant is
not as ``regular" as the preceding one, and is therefore harder to
evaluate.

\proclaim{Lemma~\TL}We have
$$L(H^-)=\det_{1\le i,j\le 2n-1}\(\frac {(2n+m-i-1)!}
{(m+i-j)!\,(2n-2i+j)!}
\left\{\matrix (n+m-\frac {j} {2})&\text {if }i\ne n\\
j&\text {if }i=n\endmatrix \right\}\).\tag\DD$$
\endproclaim
\demo{Proof}
Using the correspondence with nonintersecting lattice paths,
the number of lozenge tilings of the lower ``half hexagon" $H^-$ (see Figure
3.3) equals the
number of families $(P_1,P_2,\dots,P_{2n-1})$
of nonintersecting lattice paths (consisting of horizontal positive unit
steps and vertical negative unit steps)
in which $P_i$ runs from $A_i=(2i,i+m)$ to $E_i=(2n+i,i)$,
$i=1,2,\dots,n-1,n+1,\dots,2n-1$, whereas $P_n$ runs from
$A_n=(2n+1,n+m)$ to $E_n=(3n,n)$ (i.e., it is the starting point
$A_n$ which deviates slightly from the ``general" rule). In this
count, horizontal steps originating from any $A_i$, $i\ne n$, count
with weight $1/2$.
Again, by the main theorem of nonintersecting lattice paths,
this number is given by
the determinant
$$\det_{1\le i,j\le 2n-2}\big(\vert\Cal P(A_i\to E_j)\vert_w\big),$$
where $\vert.\vert_w$ denotes this weighted count. It is not difficult to see
that the weighted count $\vert\Cal P(A_i\to E_j)\vert_w$ of paths from
$A_i$ to $E_j$ equals $\frac {(n+m-j/2)\,(2n+m-i-1)!}
{(m+i-j)!\,(2n-2i+j)!}$ if $i\ne n$ and\linebreak $\binom
{n+m-1}{m+n-j}=\frac {j\,(n+m-1)!}
{(m+n-j)!\,j!}$ otherwise. This
establishes the Lemma.\quad \quad \qed

\enddemo

\subhead 5. Determinant evaluations\endsubhead
Here we evaluate the determinants which appear in
Lemmas~\TK, \TKa\ and \TL.

In order to compute the determinants in (\DA) and (\DAa), we utilize the
following determinant lemma from \cite{\KratAM, Lemma~2.2}.
\proclaim{Lemma~\TM}Let $X_1,\dots,X_N$, $A_2,\dots,A_N$, and $B_2,
\dots B_N$ be indeterminates. Then there holds

\vbox{
$$\multline \det_{1\le i,j\le N}\Big((X_i+A_N)\cdots(X_i+A_{j+1})
(X_i+B_j)\cdots (X_i+B_2)\Big)\\
\hskip2cm =\prod _{1\le i<j\le N} ^{}(X_i-X_j)\prod _{2\le i\le j\le N}
^{}(B_i-A_j).
\endmultline\tag\EB$$
\line{\hfil\hbox{\qed}\quad \quad }
}
\endproclaim
Now we can state, and prove, the evaluations of the determinants in (\DA)
and (\DAa), in a unified fashion.
\proclaim{Lemma~\TN}For any positive integer $n$ there holds
$$
\det_{1\le i,j\le N}\(\binom{N+m-i+1}{m+i-j}\)
=\prod _{i=1} ^{N}\frac {(N+m-i+1)!\,(i-1)!\,(2m+i+1)_{i-1}}
{(m+i-1)!\,(2N-2i+1)!}.
\tag\EC$$
\endproclaim
\demo{Proof}We take ${(N+m-i+1)!}/\big({(m+i-1)!\,(2N-2i+1)!}\big)$ out of the
$i$-th row of the determinant in (\EC), $i=1,2,\dots,N$. Thus we
obtain
$$\multline
\prod _{i=1} ^{N}\frac {(N+m-i+1)!}
{(m+i-1)!\,(2N-2i+1)!}\\
\times \det_{1\le i,j\le N}\big((m+i-j+1)\cdots (m+i-1)(N-2i+j+2)\cdots
(2N-2i+1)\big).
\endmultline$$
Taking $(-2)^{N-j}$ out of the $j$-th column, $j=1,2,\dots,N$,
we may write this as
$$\multline
(-2)^{\binom {N}2}\prod _{i=1} ^{N}\frac {(N+m-i+1)!}
{(m+i-1)!\,(2N-2i+1)!}\\
\times \det_{1\le i,j\le N}\(\Big(i-\frac {2N+1} {2}\Big)\cdots \Big(i-\frac
{N+j+2} {2}\Big)(i+m-j+1)\cdots (i+m-1)\).
\endmultline$$
Now Lemma~\TM\ can be applied with $X_i=i$, $A_j=-(N+j+1)/2$,
$B_j=m-j+1$. After some simplification one arrives at the right-hand
side of (\EC).\quad \quad \qed

\enddemo

The determinant in (\DD) evaluates as follows.
\proclaim{Lemma~\TO}For any positive integer $n$ there holds
$$\multline
\det_{1\le i,j\le 2n-1}\(\frac {(2n+m-i-1)!}
{(m+i-j)!\,(2n-2i+j)!}
\left\{\matrix (n+m-\frac {j} {2})&\text {if }i\ne n\\
j&\text {if }i=n\endmatrix \right\}\)\\
=\frac {1} {2^{3n-3}\,(n-1)!}
\prod _{i=1} ^{n}(2i-1)!^2
\prod _{i=1} ^{2n-1}\frac {(2n+m-i-1)!\,} {(m+i-1)!\,(4n-2i-1)!}
\prod _{i=1} ^{2n-2}(2m+i+1)_i\\
\times \sum _{i=0} ^{n-1}\frac {(-1)^{n-i-1}} {(2n-2i-1)}\frac
{(m+n-i)_{2i}} {i!^2}.
\endmultline\tag\EE$$
\endproclaim

\demo{Proof} The method that we use for this proof is also applied
successfully in \cite{\KratBD, \KratBG, \KratBH, \KratBI, \KrZeAA}
(see in particular the tutorial description in \cite{\KratBI, Sec.~2}).

First of all, as in the proof of Lemma~\TN, we take appropriate
factors out of the determinant. To be precise, we take
${(2n+m-i-1)!} /\big({(m+i-1)!\,(4n-2i-1)!}\big) $ out of the
$i$-th row of the determinant in (\EE), $i=1,2,\dots,2n-1$. Thus we
obtain
$$\multline
\prod _{i=1} ^{2n-2}\frac {(2n+m-i-1)!}
{(m+i-1)!\,(4n-2i-1)!}\\
\times \det_{1\le i,j\le 2n-1}\bigg((m+i-j+1)\cdots (m+i-1)(2n-2i+j+1)\cdots
(4n-2i-1)\\
\cdot\left\{\matrix (n+m-\frac {j} {2})&\text {if }i\ne n\\
j&\text {if }i=n\endmatrix\right\}\bigg).\\
\endmultline\tag\EF$$
Let us denote the determinant in (\EF) by $D(m;n)$. Using the
notation of shifted factorials, this means that
$$D(m;n):=\det_{1\le i,j\le 2n-1}\bigg((m+i-j+1)_{j-1}\,(2n-2i+j+1)_{2n-j-1}
\left\{\matrix (n+m-\frac {j} {2})&\text {if }i\ne n\\
j&\text {if }i=n\endmatrix\right\}\bigg).
\tag\EG$$
Comparison of (\EE) and (\EF) yields that (\EE) will be proved once
we are able to establish the determinant evaluation
$$\multline D(m;n)=2^{2(n-1)(n-2)}\frac {\prod _{i=1} ^{n}(2i-1)!^2}
{(n-1)!} \prod _{i=1} ^{n-1}\big((m+i)_{2n-2i}\(m+i+
{1}/{2}\)_{n-1}\big)\\
\times\sum _{i=0} ^{n-1}\frac {(-1)^{n-i-1}} {(2n-2i-1)}\frac
{(m+n-i)_{2i}} {i!^2}.
\endmultline\tag\EH$$

For the proof of (\EH) we proceed in several steps. An outline is as
follows. In the first step we show that $\prod _{i=1} ^{n-1}(m+i)_{2n-2i}$
is a factor of $D(m;n)$ 
as a polynomial in $m$. In the second step we
show that $\prod _{i=1} ^{n-1}(m+i+1/2)_{n-1}$
is a factor of $D(m;n)$.
In the third step we determine the maximal degree of $D(m;n)$ as a
polynomial in $m$, which turns out to be $(2n+1)(n-1)$. From a
combination of these three steps we are forced to conclude that
$$D(m;n)=\prod _{i=1} ^{n-1}\big((m+i)_{2n-2i}\(m+i+
{1}/{2}\)_{n-1}\big)P(m;n),\tag\EI$$
where $P(m;n)$ is a polynomial in $m$ of degree at most $2n-2$.
Then, in the fourth step we show that $P(m;n)=P(1-2n-m;n)$. And,
in the fifth step, we evaluate $P(m;n)$ at $m=0,-1,\dots,-n+1$.
Namely, for $m=0,-1,\dots,-n+1$ we show that
$$P(m;n)={{\left( -1 \right) }^n}\,{2^
       {2 \left( n -2\right)  \left(  n-1 \right) }}
\frac {  \prod_{i = 1}^{n}{{\left( 2 i-1 \right) !}^2} }
   {\left( n -1\right) !}
{{
     ({ \textstyle {1\over 2} + n}) _{m+n -1 } 
    }\over
   {2  ({ \textstyle {1\over 2} - n}) _{m+n }
}}.
\tag\EIa$$
Clearly the latter two properties determine a polynomial of maximal
degree $2n-2$ uniquely. As is easy to check, the sum in (\EH) has the
first property, too, namely that it is invariant under replacement of $m$
by $1-2n-m$. Since
in the sixth step we prove that for $m=0,-1,\dots,-n+1$ we also have
$$\sum _{i=0} ^{n-1}\frac {(-1)^{n-i-1}} {(2n-2i-1)}\frac
{(m+n-i)_{2i}} {i!^2}=(-1)^n
{{
     ({ \textstyle {1\over 2} + n}) _{m+n -1 } 
    }\over
   {2  ({ \textstyle {1\over 2} - n}) _{m+n }
}}$$
we are forced to conclude that
$$P(m;n)={2^
       {2 \left( n -2\right)  \left(  n-1 \right) }}
\frac {  \prod_{i = 1}^{n}{{\left( 2 i-1 \right) !}^2} }
   {\left( n -1\right) !}
\sum _{i=0} ^{n-1}\frac {(-1)^{n-i-1}} {(2n-2i-1)}\frac
{(m+n-i)_{2i}} {i!^2}.
\tag\EJ$$
This would finish the proof of the Lemma since a combination of (\EI)
and (\EJ) gives (\EH), and thus (\EE), as we already noted.

\bigskip
{\it Step 1. $\prod _{i=1} ^{n-1}(m+i)_{2n-2i}$ is a factor of
$D(m;n)$}. For $i$ between $1$ and $n-1$
let us consider row $2n-i$ of the determinant $D(m;n)$.
Recalling the definition (\EG) of $D(m;n)$, we see that
the $j$-th entry in this row has the form
$$(m+2n-i-j+1)_{j-1}\,(-2n+2i+j+1)_{2n-j-1}
\Big(n+m-\frac {j} {2}\Big).$$
Since $(-2n+2i+j+1)_{2n-j-1}=0$ for $j=1,2,\dots,2n-2i-1$, the first
$2n-2i-1$ entries in this row vanish. Therefore $(m+i)_{2n-2i}$ is a
factor of each entry in row $2n-i$, $i=1,2,\dots,n-1$. Hence, the
complete product $\prod _{i=1} ^{n-1}(m+i)_{2n-2i}$ divides $D(m;n)$.

\smallskip
{\it Step 2. $\prod _{i=1} ^{n-1}(m+i+1/2)_{n-1}$ is a factor of
$D(m;n)$}. Let us concentrate on a typical factor $(m+j+l+1/2)$,
$1\le j\le n-1$, $0\le l\le n-2$. We claim that for each
such factor there is a linear combination of the rows that vanishes
if the factor vanishes. More precisely, we claim that for any $j,l$
with $1\le j\le n-1$, $0\le l\le n-2$ there holds
$$\multline
{{\left( -1 \right) }^{j-1}}{{
      ({ \textstyle n - j - l -{1\over 2}}) _{j}}\over
    {{4^j} ({ \textstyle n-j-l}) _{j} }}\\
\cdot(\text {column
$(2n-2j-2l-1)$ of $D(-j-l-1/2;n)$}) \\
+    \sum_{s =  2n- j - 2 l -1}^{ 2n- 2 l -1}
     {j\choose { s+j + 2 l - 2 n + 1}} \cdot(\text {column $s$ of
$D(-j-l-1/2;n)$})=0
\endmultline\tag\EKa$$
if $j+l<n$, and
$$\multline
  \sum_{s = n-l}^{ 2 j+1}
     {{ 2 j + l - n+1}\choose {s+l - n }}
\cdot(\text {column $s$ of
$D(-j-l-1/2;n)$})\\
+{{{\left( -1 \right) }^{n-l}}\over {4^{ n-l-1}} }
{{      ({ \textstyle  j + l - n+{3\over 2}}) _{ n- l -1}  }      
     \over {({ \textstyle  j + l - n+2}) _{ n- l -1} }}\hskip3cm\\
\times
\sum_{s = 1}^{ 2 j + 2 l - 2 n+3}
         {{ 2 j + 2 l - 2 n+2}\choose { s -1}}
\cdot(\text {column $s$ of
$D(-j-l-1/2;n)$})=0
\endmultline\tag\EKb$$
if $j+l\ge n$.

In order to verify (\EKa), for $j+l<n$ we have to check
$$\multline
\sum_{s =  2n- j - 2 l -1}^{ 2n- 2 l -1}
    {j\choose { s+j + 2 l - 2 n + 1}} \\
\cdot
   \Big( n - j - l  - {s\over 2} -{1\over 2}\Big) 
    ({ \textstyle  i - j - l - s+{1\over 2}}) _{ s -1}\, 
    ({ \textstyle 2n - 2 i  + s+1}) _{ 2 n - s -1} =0,
\endmultline$$
which is (\EKa) restricted to the $i$-th row, $i\ne n$ (note
that the entry in column $2n-2j-2l-1$ of $D(-j-l-1/2;n)$ vanishes
in such a row), and
$$\multline
{{\left( -1 \right) }^{j-1}}{{
      ({ \textstyle  j + l - n+{3\over 2}}) _{ 2n- j - 2 l -2}\, 
      ({ \textstyle 2n - 2 j - 2 l -1}) _{ 2 j + 2 l+1} }\over
    {{4^j} ({ \textstyle n-j-l}) _{j} }} \\
+  \sum_{s =  2n- j - 2 l -1}^{ 2n- 2 l -1}
     {j\choose { s+j + 2 l - 2 n + 1}}
      ({ \textstyle n - j - l - s+{1\over 2}}) _{ s -1}\, 
      ({ \textstyle s}) _{2 n - s}  =0,
\endmultline$$
which is (\EKa) restricted to the $n$-th row.
Equivalently, using the standard hypergeometric notation
$${}_r F_s\!\left[\matrix a_1,\dots,a_r\\ b_1,\dots,b_s\endmatrix;
z\right]=\sum _{k=0} ^{\infty}\frac {\po{a_1}{k}\cdots\po{a_r}{k}}
{k!\,\po{b_1}{k}\cdots\po{b_s}{k}} z^k\ ,$$
this means to check
$$\multline
\kern-5pt\frac {j} {2}{{   ({ \textstyle  i + l - 2 n+{3\over 2}}) _{ 2n- j - 2 l
        -2}  \,({ \textstyle 4n-2 i - j - 2 l }) _{j + 2 l} 
     }}\,\,
{} _{3} F _{2} \!\left [ \matrix { 1 + j, -{1\over 2} - i - l
        + 2 n,-j}\\ { j, -2 i - j - 2 l + 4 n}\endmatrix ; {\displaystyle
        1}\right ]=0,\\
\endmultline\tag\ELa$$
and
$$\multline
{{\left( -1 \right) }^{j-1}}{{
      ({ \textstyle  j + l - n+{3\over 2}}) _{ 2n- j - 2 l -2}\, 
      ({ \textstyle 2n - 2 j - 2 l -1}) _{ 2 j + 2 l+1} }\over
    {{4^j} ({ \textstyle n-j-l}) _{j} }} \\
+  ({ \textstyle  l - n+{3\over 2}}) _{ 2n- j - 2 l -2}\, 
  ({ \textstyle  2n- j - 2 l -1}) _{ j + 2 l+1} \,\,
{} _{2} F _{1} \!\left [ \matrix { -{1\over 2} - l + n,-j}\\ { -1 - j - 2 l
   + 2 n}\endmatrix ; {\displaystyle 1}\right ]
=0.
\endmultline\tag\ELb$$

In order to verify (\EKb), for $j+l\ge n$ we have to check
$$\multline
 \sum_{s = n-l}^{ 2 j+1}
      {{ 2 j + l - n+1}\choose {s+l - n }}
\Big( n - j - l - {s\over 2} -{1\over 2}\Big)\\
\cdot      ({ \textstyle  i - j - l - s+{1\over 2}}) _{ s -1}\, 
      ({ \textstyle 1 - 2 i + 2 n + s}) _{ 2 n - s -1}  \\
+{{{\left( -1 \right) }^{n-l}} \over {4^{ n-l-1}} }
{{      ({ \textstyle  j + l - n+{3\over 2}}) _{ n- l -1}\, 
 }\over
      {({ \textstyle  j + l - n+2}) _{ n- l -1} }}
       \sum_{s = 1}^{ 2 j + 2 l - 2 n+3}
                {{ 2 j + 2 l - 2 n+2}\choose { s -1}} \\
\cdot    \left( n - j - l - {s\over 2} -{1\over 2}\right)
          ({ \textstyle  i - j - l - s+{1\over 2}}) _{ s -1}\, 
          ({ \textstyle 1 - 2 i + 2 n + s}) _{ 2 n - s -1} ,
\endmultline$$
which is (\EKb) restricted to the $i$-th row, $i\ne n$, and
$$\multline
  \sum_{s = n-l}^{ 2 j+1}
     {{ 2 j + l - n+1}\choose {s+l - n }}
      ({ \textstyle n - j - l - s+{1\over 2}}) _{ s -1}\, 
      ({ \textstyle s}) _{2 n - s}  \\
+{{{\left( -1 \right) }^{n-l}} \over {4^{ n-l-1}} }
{{      ({ \textstyle  j + l - n+{3\over 2}}) _{ n- l -1}\, 
 }\over
    {({ \textstyle  j + l - n+2}) _{ n- l -1} }}\kern4cm\\
\times
 \sum_{s = 1}^{ 2 j + 2 l - 2 n+3}
         {{ 2 j + 2 l - 2 n+2}\choose { s -1}}
          ({ \textstyle n - j - l  - s+{1\over 2}}) _{ s -1}\, 
          ({ \textstyle s}) _{2 n - s}  ,
\endmultline$$
which is (\EKb) restricted to the $n$-th row.
Equivalently, using hypergeometric notation, this means to check
$$\multline
\frac {\left( n - 2 j - l -1 \right)} {2}  {{ 
      ({ \textstyle {1\over 2} + i - j - n}) _{n - l -1} \,
      ({ \textstyle 3n - 2 i - l + 1}) _{n+l-1} }}\\
\times
      {} _{3} F _{2} \!\left [ \matrix { 2 + 2 j + l - n, {1\over 2} - i + j
       + n, -1 - 2 j - l + n}\\ { 1 + 2 j + l - n, 1 - 2 i - l +
       3 n}\endmatrix ; {\displaystyle 1}\right ]\\
+{{{\left( -1 \right) }^{n-l}}\over {4^{n-l-1}}}
{{
      \left( n -  j -  l -1 \right) 
      ({ \textstyle {3\over 2} + j + l - n}) _{n- l -1} \,
      ({ \textstyle 2n - 2 i + 2 }) _{2 n-2} }\over
    {({ \textstyle 2 + j + l - n}) _{n- l -1} }} \\
\times
      {} _{3} F _{2} \!\left [ \matrix { 3 + 2 j + 2 l - 2 n, {3\over 2} -
       i + j + l, -2 - 2 j - 2 l + 2 n}\\ { 2 + 2 j + 2 l - 2 n, 2 -
       2 i + 2 n}\endmatrix ; {\displaystyle 1}\right ]=0
\endmultline\tag\EMa$$
and
$$\multline
   ({ \textstyle {1\over 2} - j}) _{n- l -1} \,
   ({ \textstyle n-l}) _{n+l} \,
  {} _{2} F _{1} \!\left [ \matrix { {1\over 2} + j, -1 - 2 j - l + n}\\ { -l
    + n}\endmatrix ; {\displaystyle 1}\right ] 
\\
+{{{\left( -1 \right) }^{n-l}} \over {4^{n-l-1}}}
{{
      ({ \textstyle 1}) _{2 n-1} \,
      ({ \textstyle {3\over 2} + j + l - n}) _{n- l -1} }\over
    {({ \textstyle 2 + j + l - n}) _{n- l -1} }} \\
\times
      {} _{2} F _{1} \!\left [ \matrix { {3\over 2} + j + l - n, -2 - 2 j -
       2 l + 2 n}\\ { 1}\endmatrix ; {\displaystyle 1}\right ]=0.
\endmultline\tag\EMb$$

We start with the proof of (\ELa). We apply the contiguous relation
$$
{} _{3} F _{2} \!\left [ \matrix { a, A_1,A_2}\\ { B_1,B_2}\endmatrix ; {\displaystyle
   z}\right ]  = {} _{3} F _{2} \!\left [ \matrix { a - 1,
A_1,A_2}\\ {
    B_1,B_2}\endmatrix ; {\displaystyle z}\right ]  +
   {z    }
   {{A_1A_2 }\over    {B_1B_2}}
   {} _{3} F _{2} \!\left [ \matrix { a, A_1+1,A_2+1}\\ {B_1+1,B_2+1}\endmatrix ;
        {\displaystyle z}\right ]
$$
to the $_3F_2$-series in (\ELa). Since in this case $a-1=B_1$, parameters cancel
inside the two $_3F_2$-series on
the right hand side of the contiguous relation, leaving two $_2F_1$-series
instead. Thus, (\ELa) is turned into
$$\multline
\frac {j} {2}{{   ({ \textstyle  i + l - 2 n+{3\over 2}}) _{ 2n- j - 2 l
        -2}  \,({ \textstyle 4n-2 i - j - 2 l }) _{j + 2 l} 
     }}\bigg({} _{2} F _{1} \!\left [ \matrix {  2n - i - l -{1\over
2}},-j\\
          { 4n-2 i - j - 2 l }\endmatrix ; {\displaystyle 1}\right ]\\
 -{{\left( 2n - i - l -{1\over 2} \right) 
               }\over {4n-2 i - j - 2 l }}
{} _{2} F _{1} \!\left [ \matrix {  {1\over 2} - i - l +
               2 n,1-j}\\ { 1 - 2 i - j - 2 l + 4 n}\endmatrix ;
               {\displaystyle 1}\right ]
          \bigg)=0.
\endmultline$$
Each of the two $_2F_1$-series can be evaluated by means of the
Chu--Vandermonde summation (see \cite{\SlatAC, (1.7.7); Appendix (III.4)}),
$$
{} _{2} F _{1} \!\left [ \matrix { a, -N}\\ { c}\endmatrix ; {\displaystyle
   1}\right ]  = {{({ \textstyle c-a}) _{N} }\over
    {({ \textstyle c}) _{N} }},
\tag\EN$$
where $N$ is a nonnegative integer, and thus (\ELa) follows upon
minor simplification (the terms in big parentheses cancel each
other).

To the $_2F_1$-series in (\ELb)
Chu--Vandermonde summation (\EN) can be applied directly,
and it yields the desired result.

The verifications of (\EMa) and (\EMb) are similar. The reader will
have no difficulties to fill in the details.

This finishes the proof that the product
$\prod _{i=1} ^{n-1}(m+i+1/2)_{n-1}$ divides $D(m;n)$.

\smallskip
{\it Step 3. $D(m;n)$ is a polynomial in $m$ of maximal degree
$(2n+1)(n-1)$}. Obviously, the degree in $m$ of the $(i,j)$-entry in
the determinant $D(m;n)$ is $j$ for $i\ne n$, while it is $j-1$ for
$i=n$. Hence, in the defining expansion of the determinant, each term
has degree $\left(\sum _{j=1} ^{2n-1}j\right)-1=\binom {2n}2-1=(2n+1)(n-1)$.

\smallskip
{\it Step 4. $P(m;n)=P(1-2n-m;n)$}.
We claim that there holds the relation
$$D(m;n)=(-1)^{n-1}D(1-2n-m;n).\tag\EO$$
It is clear from the definition (\EI) of $P(m;n)$ that (\EO)
immediately implies the desired relation $P(m;n)=P(1-2n-m;n)$.

We prove (\EO) by, up to sign, transforming
the determinant $D(m;n)$ into the
determinant $D(1-2n-m;n)$ by a sequence of elementary column
operations (which, of course, leave the value of the determinant invariant).
To be precise, for $j=2n-1,2n-2,\dots,2$, we add
$$\sum _{k=1} ^{j-1}\binom {j-1}{k-1}\cdot(\text {column $k$ of
$D(m;n)$})$$
to column $j$. Thus, in the new determinant, $D_1(m;n)$ say, the
$(i,j)$-entry is
$$
\sum_{k = 1}^{j}
   {{ j -1}\choose { k -1}} \Big(  m + n-{{k}\over 2} \Big)
({ \textstyle m + i - k + 1}) _{ k -1}\, 
    ({ \textstyle 2n - 2 i + k + 1}) _{ 2n- k -1}
$$
for $i\ne n$, and
$$
\sum_{k = 1}^{j}{{ j -1}\choose { k -1}}
    ({ \textstyle m + n - k + 1}) _{ k -1} \,
    ({ \textstyle k}) _{2n- k}
$$
for $i=n$. Using hypergeometric notation, the $(i,j)$-entry of
$D_1(m;n)$ is
$$
\Big(m+n-\frac {1} {2}\Big)
{{       ({ \textstyle 2 - 2 i + 2 n}) _{ 2 n -2} }}\,\,
{} _{3} F _{2} \!\left [ \matrix { 2 - 2 m - 2 n, 1 - i - m, 1 - j}\\ {
      1 - 2 m - 2 n, 2 - 2 i + 2 n}\endmatrix ; {\displaystyle 1}\right ]
\tag\EPa$$
for $i\ne n$, and
$$
  ({ \textstyle 1}) _{ 2 n -1} \,\,
{} _{2} F _{1} \!\left [ \matrix { 1 - m - n,1-j}\\ { 1}\endmatrix ;
   {\displaystyle 1}\right ]
\tag\EPb$$
for $i=n$. The $_3F_2$-series in (\EPa) can be evaluated in the same
way as we evaluated the $_3F_2$-series in (\ELa). The $_2F_1$-series
in (\EPb) is easily evaluated by Chu--Vandermonde summation (\EN).
Thus, we obtain that the $(i,j)$-entry of $D_1(m;n)$ is given by
$$
(-1)^j\Big( 1 - {j\over 2} - m - n \Big) 
  ({ \textstyle 2 + i - j - m - 2 n}) _{ j -1}\, 
  ({ \textstyle 1 - 2 i + j + 2 n}) _{ 2n- j -1}
\tag\EQa$$
for $i\ne n$, and by
$$
{{\left( -1 \right) }^{ j -1}} j\, 
  ({ \textstyle 2 - j - m - n}) _{ j -1} \,
({ \textstyle j+1}) _{2n-j -1}
\tag\EQb$$
for $i=n$. Now, the expression (\EQa) is exactly $(-1)^j$ times the
$(i,j)$-entry of $D(1-2n-m;n)$, while the expression (\EQb) is
exactly $(-1)^{j-1}$ times the $(n,j)$-entry of $D(1-2n-m;n)$. Hence,
relation (\EO) follows immediately, implying $P(m;n)=P(1-2n-m;n)$, as
we already noted.

\smallskip
{\it Step 5. Evaluation of $P(m;n)$ at $m=0,-1,\dots,-n+1$}.
The polynomial $P(m;n)$ is defined by means of (\EI),
$$D(m;n)=\prod _{i=1} ^{n-1}\big((m+i)_{2n-2i}\(m+i+
{1}/{2}\)_{n-1}\big)P(m;n).\tag\EIb$$
So, what we
would like to do is to set $m=-e$, $e$ being one of
$0,1,\dots,n-1$, evaluate $D(-e;n)$, divide both sides of (\EIb) by
the product on the right-hand side of (\EIb), and get the evaluation
of $P(m;n)$ at $m=-e$. However, the product on the right-hand side of
(\EIb) unfortunately (usually) {\it is zero} for $m=-e$, $0\le e\le
n-1$. Therefore we have to find a way around this difficulty.

Fix an $e$ with $0\le e\le n-1$. Before setting $m=-e$ in (\EIb), we
have to cancel $(m+e)^e$ on the right-hand side of (\EIb). To accomplish this,
we have to ``generate" these factors on the left-hand side. We do
this by adding
$$\sum _{k=j+1} ^{2e+2j-2n}\binom {2e+j-2n}{k-j}\cdot(\text {column $k$ of
$D(m;n)$})$$
to column $j$, $j=2n-2e,2n-2e+1,\dots,2n-e-1$.
Thus, in the new determinant the entry in the $i$-th
row in such a column is
$$
\sum_{k = j}^{2e+2j-2n}\Big( m + n - {k\over 2} \Big) 
    {{2 e + j - 2 n}\choose {k-j}}
    ({ \textstyle m + i - k+1}) _{ k -1}\, 
    ({ \textstyle 2n - 2 i + k+1}) _{ 2 n - k -1}
$$
if $i\ne n$, and
$$
\sum_{k = j}^{2e+2j-2n}{{2 e + j - 2 n}\choose {k-j }}
    ({ \textstyle  m + n - k+1}) _{ k -1}\,  ({ \textstyle k}) _{2 n - k}
$$
if $i=n$. In hypergeometric terms this is
$$\multline
\Big(m+n-\frac {j} {2}\Big)
{{       ({ \textstyle m + i - j + 1}) _{ j -1}\, 
     ({ \textstyle 2n - 2 i + j + 1}) _{ 2n- j -1} }}\\
\times
{} _{3} F _{2} \!\left [ \matrix { 1 + j - 2 m - 2 n, -i + j - m, -2 e
      - j + 2 n}\\ { j - 2 m - 2 n, 1 - 2 i + j + 2 n}\endmatrix ;
      {\displaystyle 1}\right ]
\endmultline$$
if $i\ne n$, and
$$
  ({ \textstyle m+n - j + 1}) _{ j -1}\,
({ \textstyle j}) _{2n-j}
\,{} _{2} F _{1} \!\left [ \matrix {  j - m - n,-2 e - j + 2 n}\\ {
   j}\endmatrix ; {\displaystyle 1}\right ] 
$$
if $i=n$.
Again, the $_3F_2$-series can be evaluated in the same way as before
the $_3F_2$-series in (\ELa), while the $_2F_1$-series is evaluated
by means of the Chu--Vandermonde summation (\EN). Thus, we obtain
that the $(i,j)$-entry, $2n-2e\le j\le 2n-e-1$,
of the modified determinant is given by
$$
\left(  m+e \right)    ({ \textstyle m + i - j + 1}) _{ j -1}\, 
  ({ \textstyle m + 2 n-i}) _{2 e + j - 2 n} \,
({ \textstyle 2 e - 2 i + 2 j+1}) _{ 4n- 2 e - 2 j -1}
$$
if $i\ne n$, and by
$$
\left(  m+e \right)  \, 
  ({ \textstyle m+n - j + 1}) _{ e + j - n -1} \,
({ \textstyle m+e+1}) _{ e + j - n -1}\,
  ({ \textstyle 2 e + 2 j - 2 n}) _{4n-2 e - 2 j } 
$$
if $i=n$. Clearly, $(m+e)$ is a factor of each entry in the $j$-th
column of the modified determinant, $2n-2e\le j\le 2n-e-1$.
Therefore, we may take $(m+e)$ out of the $j$-th column,
$j=2n-2e,2n-2e+1,\dots,2n-e-1$. The remaining determinant, $D_2(m;n)$
say, is then defined as $D_2(m;n)=\det_{1\le i,j\le 2n-1}(E_{ij})$,
where for $1\le j\le 2n-2e-1$ and for $2n-e\le j\le 2n-1$ the entry
$E_{ij}$ is
given by
$$E_{ij}=(m+i-j+1)_{j-1}\,(2n-2i+j+1)_{2n-j-1}
\left\{\matrix (n+m-\frac {j} {2})&\text {if }i\ne n\\
j&\text {if }i=n\endmatrix\right.\tag\ERa$$
(i.e., $E_{ij}$ equals the $(i,j)$-entry of $D(m;n)$, as given by
(\EG), in that case), and for $2n-2e\le j\le 2n-e-1$ the entry $E_{ij}$ is
given by
$$E_{ij}=\cases  ({ \textstyle m + i - j + 1}) _{ j -1}\, 
  ({ \textstyle m + 2 n-i}) _{2 e + j - 2 n} \,
({ \textstyle 2 e - 2 i + 2 j+1}) _{ 4n- 2 e - 2 j -1} &\text {if }i\ne n\\
  ({ \textstyle m+n - j + 1}) _{ e + j - n -1} \,
({ \textstyle e + m+1}) _{ e + j - n -1}\, 
  ({ \textstyle 2 e + 2 j - 2 n}) _{4n-2 e - 2 j } 
&\text {if }i=n.
\endcases\tag\ERb$$
Due to the manipulations that we did, the new determinant $D_2(m;n)$
is related to the original determinant $D(m;n)$ by
$$D(m;n)=(m+e)^e\,D_2(m;n).$$
Substituting this into (\EIb), and rearranging terms, we get
$$\multline
P(m;n)=
D_2(m;n)\prod _{i=1} ^{e}\big((m+i)_{e-i}\,(m+e+1)_{2n-i-e-1}\(m+i+
{1}/{2}\)_{n-1}\big)^{-1}\\
\times \prod _{i=e+1} ^{n-1}\big((m+i)_{2n-2i}\(m+i+
{1}/{2}\)_{n-1}\big)^{-1}.
\endmultline\tag\ES$$
Now we may safely set $m=-e$. So, what we need in order to obtain the
evaluation of $P(m;n)$ at $m=-e$ is the evaluation of the determinant
$D_2(-e;n)$.

In order to determine the evaluation of $D_2(-e;n)$, we observe that
$D_2(-e;n)$ has a block form which is sketched in Figure~\FZ.

\vskip10pt
\vbox{\noindent
$$
\PfadDicke{.5pt}
\thinlines
\Pfad(0,0),111111\endPfad
\Pfad(0,2),111111\endPfad
\Pfad(0,4),111111\endPfad
\Pfad(0,6),111111\endPfad
\Pfad(0,0),222222\endPfad
\Pfad(2,0),222222\endPfad
\Pfad(4,0),222222\endPfad
\Pfad(6,0),222222\endPfad
\Pfad(0,4),4444\endPfad
\Label\l{\kern18pt$\fourteenpoint$0}(1,1)
\Label\l{\raise10pt\hbox{$\kern30pt$\fourteenpoint$*$}}(3,1)
\Label\l{\kern18pt$\fourteenpoint$0}(5,3)
\Label\l{\kern18pt$\fourteenpoint$0}(5,1)
\Label\l{\raise10pt\hbox{$\kern30pt$\fourteenpoint$0$}}(1,3)
\Label\l{\raise-15pt\hbox{$\kern5pt$\fourteenpoint$*$}}(1,3)
\Label\l{\kern18pt$\fourteenpoint$*}(1,5)
\Label\l{\raise-15pt\hbox{$\kern5pt$\fourteenpoint$0$}}(3,1)
\Label\l{\kern18pt$\fourteenpoint$*}(3,3)
\Label\l{\kern18pt$\fourteenpoint$*}(3,5)
\Label\l{\kern18pt$\fourteenpoint$*}(5,5)
\Label\r{i=1\hphantom{e+{}}}(7,6)
\Label\r{i=e+1}(7,4)
\Label\r{\kern5pt i=2n-e}(7,2)
\Label\r{\kern5pt i=2n\hphantom{{}-1}}(7,0)
\Label\o{j=1}(0,6)
\Label\u{j=}(2,0)
\Label\u{\raise15pt\hbox{$2n-2e$}}(2,-1)
\Label\o{\raise-15pt\hbox{$j=$}}(4,7)
\Label\o{2n-e}(4,6)
\Label\u{j=}(6,0)
\Label\u{\raise15pt\hbox{$2n$}}(6,-1)
\hskip3cm
$$
\centerline{\eightpoint The block form of $D_2(-e;n)$}
\vskip8pt
\centerline{\eightpoint Figure \FZ}
}
\vskip10pt
\noindent
The figure has to be read according to the following convention:
If a block is
bounded by horizontal lines marked as $i=h_1$ and $i=h_2$ and
vertical lines marked as $j=v_1$ and $j=v_2$, then the block consists
of the entries that are in rows $i=h_1,h_1+1,\dots,h_2-1$ and
columns $j=v_1,v_1+1,\dots,v_2-1$. It is an easy task to check from
the definitions (\ERa) and (\ERb) of the entries of $D_2(m;n)$ that
indeed the lower left block of $D_2(-e;n)$, consisting of the entries in
rows $i=2n-e,2n-e+1,\dots,2n-1$ and columns
$j=1,2,\dots,2n-2e-1$, and the lower right block, consisting of
the entries in rows $i=e+1,e+2,\dots,2n-1$ and
columns $j=2n-e,2n-e+1,\dots,2n-1$ are blocks of zeroes. Hence,
the determinant $D_2(-e;n)$ factors into the product
$$\det(B_1)\,\det(B_2)\,\det(B_3),$$
where $B_1$ is the middle left block, consisting of the entries in
rows $i=e+1,e+2,\dots,2n-e-1$ and
columns $j=1,2,\dots,2n-2e-1$,
where $B_2$ is the lower middle block, consisting of the entries in
rows $i=2n-e,2n-e+1,\dots,2n-1$ and
columns $j=2n-2e,2n-2e+1,\dots,2n-e-1$, and
where $B_3$ is the upper right block, consisting of the entries in
rows $i=1,2,\dots,e$ and
columns $j=2n-e,2n-e+1,\dots,2n-1$.

As is indicated in Figure~\FZ, the blocks $B_1$ and $B_2$ are lower
and upper triangular matrices, respectively. Hence, their
determinants are easily computed.

The determinant of the block $B_3$ is the
determinant
$$\det_{1\le i,j\le e}\big((-2n+i-j+2)_{2n-e+j-2}\,(4n-e-2i+j)_{e-j}
\Big(\frac {1-e-j} {2}\Big)\big).$$
(All the entries of $B_3$ are given by the first case of
formula (\ERa) because $e\le n-1$.) We take $(-2n+i+1)_{2n-e-1}$
out of the $i$-th row of this determinant, $i=1,2,\dots,e$,
and $(-2)^{e-j}\(\frac {1-e-j} {2}\)$
out of the $j$-th column, $j=1,2,\dots,e$. Thus we obtain
$$\multline
(-2)^{\binom e2}\prod _{i=1} ^{e}\Big(\frac {1-e-i} {2}\Big)(-2n+i+1)_{2n-e-1}\\
\times
\det_{1\le i,j\le e}\(
\Big(i-2n+\frac {1} {2}\Big)\cdots \Big(i-2n+\frac {e-j} {2}\Big)
(i-2n-j+2)\cdots (i-2n)
\).
\endmultline$$
The latter determinant is easily evaluated using Lemma~\TM\ with
$N=e$, $X_i=i$, $A_j=-2n+(e-j+1)/2$, $B_j=-2n-j+2$.

This finishes the desired evaluation of $D_2(-e;n)$, and, via (\ES),
of $P(-e;n)$ for $e=0,1,\dots,n-1$.
If everything is put together and simplified, the
result is exactly (\EIa).

\smallskip
{\it Step 6. Evaluation of the sum on the right-hand side of\/ {\rm
(\EE)} at $m=0,-1,\dots,-n+1$}.
We claim that if $0\le e\le n-1$ we have
$$
\sum _{i=0} ^{n-1}\frac {(-1)^{n-i-1}} {(2n-2i-1)}\frac
{(n-e-i)_{2i}} {i!^2}=
{{\left( -1 \right) }^n}
{{ ({ \textstyle {1\over 2} + n}) _{n- e -1}}\over
     {2\, ({ \textstyle {1\over 2} - n}) _{n-e} }}.
\tag\ET$$
This is seen by first rewriting the sum in (\ET) as
$${{{{\left( -1 \right) }^n}  }\over {(1 - 2 n)}}
 \sum_{i = 0}^{n-1}
        {{({ \textstyle {1\over 2} - n}) _{i} \,
            ({ \textstyle 1 + e - n}) _{i} \,({ \textstyle n-e}) _{i} }
           \over {{{({ \textstyle 1}) _{i} }^2}\,
            ({ \textstyle {3\over 2} - n}) _{i} }} .
$$
Since $0\le e\le n-1$,
the term $(1+e-n)_i$ will make the sum terminate at $i=n-e-1$, and so
at $i=n-1$ latest. Therefore we may extend the range of summation to
{\it all\/} nonnegative numbers $i$, without altering the sum. Then
we can write the sum in hypergeometric notation as
$${{ {{\left( -1 \right) }^n} }\over {(1-2n)}}
{} _{3} F _{2} \!\left [ \matrix { {1\over 2} - n,
     n -e, 1 + e - n}\\ { 1, {3\over 2} - n}\endmatrix ; {\displaystyle
      1}\right ].
$$
This $_3F_2$-series can be evaluated by means of the
Pfaff--Saalsch\"utz summation (see \cite{\SlatAC, (2.3.1.3); Appendix
(III.2)}),
$$
{} _{3} F _{2} \!\left [ \matrix { a, b, -N}\\ { c, 1 + a + b - c -
   N}\endmatrix ; {\displaystyle 1}\right ]  =
  {{({ \textstyle c-a}) _{N}  ({ \textstyle c-b}) _{N} }\over
    {({ \textstyle c}) _{N}  ({ \textstyle c-a-b}) _{N} }},
$$
where $N$ is a nonnegative integer. Thus we arrive at the right-hand
side of (\ET).

\smallskip
This completes the proof of the Lemma.\quad \quad \qed
\enddemo

\enddocument